\documentclass[10pt,reqno]{amsart}
\usepackage{amsmath}
\usepackage{amssymb}
\usepackage{amsthm}

\newcommand{\nabp}{\nab \hspace{-1.05ex}
\rule[.5ex]{.2ex}{.8ex}   \hspace{1.05ex}}

\newcommand{\spat}{\hspace{4ex}}

\newcommand{\cA}{{\mathcal A}}

\newcommand{\impl}{\Upsilon}
\newcommand{\impla}{\Xi}
\newcommand{\flip}{\raisebox{.45ex}[0pt][0pt]{$\chi$}}

\newcommand{\Ga}{\Gamma}

\newcommand{\nab}{\nabla}

\newcommand{\cL}{{\mathcal L}}

\newcommand{\cN}{{\mathcal N}}
\newcommand{\cM}{{\mathcal M}}
\newcommand{\cP}{{\mathcal P}}

\newcommand{\cK}{{\mathcal K}}

\newcommand{\ot}{\otimes}

\newcommand{\la}{\Lambda}

\newcommand{\om}{\omega}
\newcommand{\io}{\iota}
\newcommand{\vfi}{\varphi}

\newcommand{\al}{\alpha}
\newcommand{\be}{\beta}

\newcommand{\sde}{\delta}
\newcommand{\de}{\Delta}

\newcommand{\si}{\sigma}

\newcommand{\Nfi}{{\mathcal N}_{\vfi}}

\newcommand{\lafi}{\la_\vfi}
\newcommand{\laps}{\la_\psi}

\newcommand{\pifi}{\pi_\vfi}
\newcommand{\pips}{\pi_\psi}

\newcommand{\otd}{\stackrel{.}{\ot}}

\newcommand{\C}{\mathbb C}
\newcommand{\R}{\mathbb R}

\newcommand{\N}{\mathbb N}

\newcommand{\Ext}{\hspace{-1.5ex}\raisebox{-0.5ex}[0pt][0pt]{\scriptsize\fontshape{n}\selectfont ext}}

\newcommand{\cU}{{\mathcal U}}

\newcommand{\cT}{{\mathcal T}}
\newcommand{\cJ}{{\mathcal J}}

\newcommand{\cst}{\text{C}$\hspace{0.1mm}^*$}

\newtheorem{definition}{Definition}[section]
\newtheorem{proposition}[definition]{Proposition}
\newtheorem{lemma}[definition]{Lemma}
\newtheorem{corollary}[definition]{Corollary}
\newtheorem{remark}[definition]{Remark}

\newtheorem{notation}[definition]{Notation}
\newtheorem{result}[definition]{Result}

\begin{document}

\begin{center}
\Huge\bf Induced corepresentations of locally compact quantum groups  \end{center}

\bigskip\bigskip\bigskip

\begin{center}
\rm Johan Kustermans

Department of Mathematics

KU Leuven

Celestijnenlaan 200B

3001 Heverlee

Belgium

\medskip

e-mail : johan.kustermans@wis.kuleuven.ac.be

\bigskip\bigskip

\bf March 2000 \rm
\end{center}

\bigskip\medskip

We introduce the construction of induced corepresentations in the setting of locally compact quantum groups and prove
that the resulting induced corepresentations are unitary under some mild integrability condition. We also establish a
quantum analogue of the classical bijective correspondence between quasi-invariant measures and certain measures on the
larger locally compact group.

\bigskip\medskip


\bigskip\medskip\medskip

\section*{Introduction}

\bigskip

Consider a closed subgroup $H$ of a locally compact group $G$ together with a strongly continuous unitary
representation $u$ of $H$ on a Hilbert space $K$. The construction of the induced representation uses these three
ingredients to supply a new strongly continuous unitary representation $\rho$ of the larger group $G$ on a new Hilbert
space $\cK$.

Induced group representations for finite groups where first introduced by Frobenius in 1898. The theory of induced
representations for locally compact groups was initiated by Mackey in 1949 (see \cite{Mack1}, \cite{Mack2} and
\cite{Mack3}). He restricted himself to locally compact groups satisfying the second axiom of countability and
separable Hilbert spaces. The general case was treated by Blattner in 1961 (see \cite{Blatt}).

The role of the theory of induced representations in the classical theory of locally compact groups can hardly be
underestimated. This construction is for instance one of the primary instruments to construct the different series of
special representations of various non-compact locally compact groups.

\medskip

Building on the work of  Enock \&  Schwartz and  Kac \& Vainerman in the theory of Kac algebras (see \cite{E}), Baaj \&
Skandalis, Kirchberg, Woronowicz and Masuda \& Nakagami in the more general setting of quantum groups,  S. Vaes and the
author developed a relatively simple definition of a locally compact quantum group. These locally compact quantum
groups can present themselves in different forms. Two of these forms are formulated in an \cst-algebraic setting (see
\cite{VaKust} and \cite{Kust1}) but in this paper we are mainly in the von Neumann algebraic setting (see
\cite{VaKust2}).

\medskip

Building on this definition, it is important to further develop the proper theory of locally compact quantum groups and
in this way obtain a theory that is as rich as the classical theory of locally compact groups. In this light, it is
very natural to generalize the theory of induced representations to the quantum group setting. This paper deals with
the first step in this generalization, the construction of the induced corepresentation and establishing its unitarity.

Starting from two von Neumann algebraic quantum groups $(M,\de)$ and $(N,\de)$ together with a special action $\al : M
\rightarrow M \ot N$ (that mimics the role of the obvious action of the group $H$ on $G$) and a unitary
corepresentation $U$ of the 'smaller' quantum group $(N,\de)$, we construct a new corepresentation $\rho$ of $(M,\de)$.
Under a mild integrability condition, we prove that this new corepresentation $\rho$ is again unitary. The difficulty
of proving the unitarity of the corepresentation $\rho$ is even more pronounced than  proving the unitarity of
multiplicative unitaries popping up in the axiomatic approach to quantum groups.

The further development of the theory of induced corepresentations will be postponed to later papers.

\medskip

Thanks to the theory of normal faithful semifinite weights on von Neumann algebras, the existence of the quantum
version of quasi-invariant measures in the classical theory becomes a mere observation in the quantum setting. Although
the classical quasi-invariant measures provide Radon Nikodym derivatives, the role of these Radon Nikodym derivatives
is much less important in the quantum setting because they seize to exist in general. Since these Radon Nikodym
derivatives appear prominently in the definition of the induced representations in the classical theory, another
approach had to be found.

\smallskip

In this respect, the results in \cite{Va1} come to the rescue and play a vital role in getting the theory going. In
this paper, S. Vaes proves that every action of a von Neumann algebraic quantum group on another von Neumann algebra is
implemented by a canonical unitary corepresentation (this is a generalization of a classical result of U. Haagerup).

\medskip

The paper is organized as follows. In section 1 we give a short overview of weight theory on von Neumman algebras. The
second section contains the definition of a von Neumann algebraic quantum group. We fix the data that we will need for
the rest of the paper in section 3. The carrier space of the induced corepresentation is constructed in section 4,
together with a useful characterization of the tensor product of this carrier space with any Hilbert space. The induced
corepresentation itself is constructed in section 5. Under an integrability condition, we provide an important dense
subspace of the carrier space in section 6. In section 7, we prove the unitarity of the induced corepresentation under
this integrability condition. In the last section we establish a bijective correspondence between the quasi-invariant
weights and a certain class of weights on the 'larger' von Neumann algebra $M$.

\medskip\medskip

We end this introduction with some notations and conventions. If $V$ is a normed space and $L$ is a subset of $V$, then
$\langle L \rangle$ will denote the linear span of $L$, $[L]$ will denote the closed linear span of $L$.  For any set
$I$, we define $F(I)$ to be the set of all finite subset of $I$ and we turn $F(I)$ into a directed set through the
inclusion relation. The identity map will  be denoted by $\io$. If $V$, $W$ are two vector spaces, the algebraic tensor
product is denoted by $V \odot W$ (we will also use $\odot$  for the tensor product of linear maps).

\smallskip

Consider two von Neumann algebras $M$ and $N$. The von Neumann algebraic tensor product of $M$ and $N$ will be denoted
by $M \ot N$ (again, we will also use $\ot$ for the tensor product of sufficiently nice normal linear maps). We denote
the flip automorphism from $M \ot N$ to $N \ot M$ by $\flip$.

\bigskip

\subsection*{Acknowledgment} I would like to thank the research group in Operator Algebras in Leuven for providing an
excellent working environment. Discussions with Stefaan Vaes about  quantum group actions turned out to be very
fruitful for this work and were very much appreciated by the author.

\bigskip\medskip

\section{Preliminaries on weights and operator valued weights}

\bigskip

Before starting our discussion on weights, we would first like to mention the following useful property. Consider
Hilbert spaces $H$,$K$ and elements $X,Y \in B(H \ot K)$, $v,w \in K$ and $(e_i)_{i \in I}$ an othonormal basis for
$K$. Then the net
\begin{equation}
\bigl(\, \sum_{i \in J} \, (\io \ot \om_{e_i,w})(X)\,(\io \ot \om_{v,e_i})(Y)\,\bigr)_{J \in F(I)} \label{prel.eq1}
\end{equation}
is bounded and converges strongly$^*$ to $(\io \ot \om_{v,w})(X Y)$  (see e.g. lemma 9.5 of \cite{VaKust} for a proof).

\bigskip

We assume that the reader is familiar with the theory of normal faithful semi-finite weights (in short, n.s.f. weights)
on von Neumann algebras. Nevertheless, let us fix some notations. So let $\vfi$ be a n.f.s. weight on a von Neumann
algebra. Then we define the following sets:
\begin{enumerate}
\item $\cM_\vfi^+ = \{ \, x \in M^+ \mid \vfi(x) < \infty \, \}$, so $\cM_\vfi^+$ is a hereditary cone in $M^+$,
\item $\cN_\vfi = \{ \, x \in M \mid x^*x \in \cM_\vfi^+ \, \}$, so $\cN_\vfi$ is a left ideal in $M$,
\item $\cM_\vfi =$ the linear span of $\cM_\vfi^+$ in $M$, so $\cM_\vfi$ is a sub $^*$-algebra of $M$.
\end{enumerate}
There exists a unique linear map $F : \cM_\vfi \rightarrow \C$ such that $F(x) = \vfi(x)$ for all $x \in \cM_\vfi^+$.
For all $x \in \cM_\vfi$, we set $\vfi(x) = F(x)$.

\smallskip

A GNS-construction for $\vfi$ is a triple $(H_\vfi,\pifi,\lafi)$, where $H_\vfi$ is a Hilbert space, $\pifi : M
\rightarrow B(H_\vfi)$ is a normal $^*$-homomorphism and $\lafi : \cN_\vfi \rightarrow H_\vfi$ is a $\si$-strong$^*$
closed linear map with dense range such that
\begin{enumerate}
\item $\langle \lafi(x), \lafi(y) \rangle = \vfi(y^* x)$ for all $x,y \in \cN_\vfi$,
\item $\lafi(x\, y) = \pifi(x) \lafi(y)$ for all $x \in M$ and $y \in \cN_\vfi$.
\end{enumerate}
Such a GNS-construction always exist. The modular group of $\vfi$ will be denoted by $\si^\vfi$, the modular operator
by $\nab_\vfi$ and the modular conjugation by $J_\vfi$ (these last two objects are defined with respect to the
GNS-construction). Recall that $\lafi(\Nfi^* \cap \Nfi)$ is a core for the operator $J_\vfi \nab_\vfi^{\frac{1}{2}}$
and $(J_\vfi\nab_\vfi^{\frac{1}{2}})\, \lafi(x) = \lafi(x^*)$ for all $x \in \Nfi \cap \Nfi^*$.

\bigskip

Consider two von Neumann algebras $M$, $N$. Let $\vfi$ be a n.f.s. weight on $M$ with GNS-construction
$(H_\vfi,\pifi,\lafi)$ and let $\psi$ be a n.s.f. weight on $N$ with GNS-construction $(H_\psi,\pips,\laps)$. The
tensor product weight $\vfi \ot \psi$ is a n.s.f. weight on $M \ot N$ (see e.g. definition 8.2 of \cite{Stramod} for a
definition). This tensor product weight has a GNS-construction $(H_\vfi \ot H_\psi,\pifi \ot \pips,\lafi \ot \laps)$
where $\lafi \ot \laps : \cN_{\vfi \ot\psi} \rightarrow H_\vfi \ot H_\psi$ is the $\si$-strong$^*$ closure of $\lafi
\odot \laps : \cN_{\vfi} \odot \cN_{\psi} \rightarrow H_\vfi \ot H_\psi$.

\bigskip

Let $M$ be any von Neumann algebra. For the definition of the extended positive part $M^+\Ext$ of $M$ we refer to
definition 1.1 of \cite{Haa1}. For $T \in M^+\Ext$ and $\om \in M_*^+$, we set $\om(T) = T(\om) \in [0,\infty]$. Recall
that there exists an embedding $M^+  \hookrightarrow M^+\Ext : x \mapsto x^\sharp$ such that $x^\sharp(\om) = \om(x)$
for all $x \in M^+$ and $\om \in M_*^+$. We will use this embedding to identify $M^+$ as a subset of $M^+\Ext$.

\smallskip

Consider a von Neumann algebra $M$ and a sub von Neumann algebra $N$ of $M$. The definition of an operator valued
weight from $M$ to $N$ is given in definition 2.1 of \cite{Haa1}.

\medskip

Now consider two von Neumann algebras $M$ and $N$ and a n.f.s. weight $\vfi$ on $M$ with GNS-construction
$(H_\vfi,\pifi,\lafi)$. We identify $N$ with $1 \ot N$ as a sub von Neumann algebra of $M \ot N$ to get into the
framework of operator valued weights. The operator valued weight $\vfi \ot \io : (M \ot N)^+ \rightarrow N^+\Ext$ is
defined in such a way that for $x \in (M \ot N)^+$, we have that $$\om\bigl((\vfi \ot \io)(x)\bigr) = \vfi\bigl((\io
\ot \om)(x)\bigr)\ .$$ As for weights we define the following sets:
\begin{enumerate}
\item $\cM_{\vfi \ot \io}^+ = \{ \, x \in (M \ot N)^+ \mid (\vfi \ot \io)(x) \in N^+ \, \}$, so $\cM_{\vfi \ot \io}^+$ is a hereditary cone of $(M \ot N)^+$,
\item $\cN_{\vfi \ot \io} = \{ \, x \in M \ot N \mid x^*x \in \cM_{\vfi \ot \io}^+ \, \}$, so $\cN_{\vfi \ot \io}$ is a left ideal in $M \ot N$,
\item $\cM_{\vfi \ot \io} =$ the linear span of $\cM_{\vfi \ot \io}^+$ in $M \ot N$, so $\cM_{\vfi \ot \io}$ is a sub $^*$-algebra of $M \ot N$.
\end{enumerate}
There exists a unique linear map $G : \cM_{\vfi \ot \io} \rightarrow N$ such that $G(x) = (\vfi \ot \io)(x)$ for all $x
\in \cM_{\vfi \ot \io}^+$. For all $x \in \cM_{\vfi \ot \io}$, we set $(\vfi \ot \io)(x) = G(x)$. Let $a \in \cM_\vfi$
and $b \in N$. Then it is easy to see that $a \ot b$ belongs to $\cM_{\vfi \ot \io}$ and $(\vfi \ot \io)(a \ot b) =
\vfi(a) \, b$.

\smallskip

Thanks to the remark after lemma 1.4 of \cite{Haa1}, we also have the following characterization of $\cM_{\vfi \ot
\io}^+$: Let $x \in (M \ot N)^+$, then $x$ belongs to $\cM_{\vfi \ot \io}^+$ $\Leftrightarrow$ $\vfi((\io \ot \om)(x))
< \infty$ for all $\om \in M^+_*$.

Let $x \in \cN_{\vfi \ot \io}$ and $\om \in N_*$. The inequality $(\io \ot \om)(x)^* (\io \ot \om)(x) \leq \|\om\| \,
(\io \ot |\om|)(x^* x)$ implies that $(\io \ot \om)(x) \in \cN_\vfi$ and $$\|\lafi((\io \ot \om)(x))\| \leq
\|\om\|\,\|(\lafi \ot \io)(x)\| \ .$$

\medskip

As for weights, there exists also a KSGNS-construction for this operator valued weight $\vfi \ot \io$. For this, let
$H$ denote the Hilbert space on which $N$ acts.

\begin{definition}
There exists a unique linear map $\lafi \ot \io : \cN_{\vfi \ot \io} \rightarrow B(H,H_\vfi \ot H)$ such that
\newline $\langle (\lafi \ot \io)(x) \, v , \lafi(a) \ot w \rangle  = \langle (\vfi \ot \io)((a^* \ot 1)x)\,v , w \rangle$
for all $x \in \cN_{\vfi \ot \io}$, $a \in \cN_\vfi$ and $v,w \in H$. Moreover, the following properties hold:
\begin{enumerate}
\item $(\lafi \ot \io)(x\,y) = (\pifi \ot \io)(x)\,(\lafi \ot \io)(y)$ for all $x \in M \ot N$ and $y \in \cN_{\vfi \ot \io}$
\item $(\lafi \ot \io)(y)^* (\lafi \ot \io)(x) = (\vfi \ot \io)(y^* x)$ for all $x,y \in \cN_{\vfi \ot \io}$.
\item $H_\vfi \ot H = [\,(\lafi \ot \io)(x)\,v \mid x \in \cN_{\vfi \ot \io}, v \in H\,]$.
\end{enumerate}
\end{definition}

We have of course that $(\lafi \ot \io)(a \ot b) \, v = \lafi(a) \ot  b\,v$ for all $a \in \cN_{\vfi}$, $b \in N$ and
$v \in H$.

\smallskip

The proof of the existence of the map $\lafi \ot \io$ and most of its properties can be extracted from  the following
result:

\begin{result} \label{prel.res1}
Consider $x \in \cN_{\vfi \ot \io}$, $v \in H$ and an orthonormal basis $(e_i)_{i \in I}$ for $H$. Then
\newline $\sum_{i \in I} \|\lafi((\io \ot \om_{v,e_i})(x))\|^2 < \infty$ and $$(\lafi \ot \io)(x)\,v = \sum_{i \in I}
\lafi((\io \ot \om_{v,e_i})(x)) \ot e_i \ .$$
\end{result}

It is also possible to show that $\lafi \ot \io : \cN_{\vfi \ot \io} \rightarrow B(H,H_\vfi \ot H)$ is closed for the
$\si$-strong$^*$ topology on $M \ot N$ and the strong topology on $B(H,H_\vfi \ot H)$ (see proposition 3.23 of
\cite{VaKust1} for a similar result in the \cst-algebra setting). One can even prove the following fact:

Let $x \in \cN_{\vfi \ot \io}$. Then their exists a net $(x_i)_{i \in I}$ in $\cN_\vfi \odot N$ such that  $\|x_i\|
\leq \|x\|$ and $\|(\lafi \ot \io)(x_i)\| \leq \|(\lafi \ot \io)(x)\|$ for all $i \in I$, $(x_i)_{i \in I}$ converges
strongly$^*$ to $x$ and $\bigl(\,(\lafi \ot \io)(x_i)\,\bigr)_{i \in I}$ converges strongly$^*$ to $(\lafi \ot \io)(x)$
(techniques for the proof of this fact can be found in the proofs of proposition 3.28 of \cite{VaKust1} and proposition
7.9 of \cite{Kust}).

\smallskip

We should also mention the following dominated convergence property (see also proposition 3.24 of \cite{VaKust1}): Let
$(x_i)_{i \in I}$ be a net in $\cM_{\vfi \ot \io}^+$ and $x$ an element in $\cM_{\vfi \ot \io}$ such that $x_i \leq x$
for all $i \in I$ and $(x_i)_{i \in I}$ converges strongly to $x$. Then $\bigl(\,(\vfi \ot \io)(x_i)\,\bigr)_{i \in I}$
converges strongly to $(\vfi \ot \io)(x)$.

\medskip

It is clear that we also define $\io \ot \vfi$ and $\io \ot \lafi$ in a similar way. But we will also use the following
variations. Let $L$,$M$,$N$ be three von Neumann algebras and let $\vfi$ be a n.s.f. weight on $L$ with
GNS-construction $(H_\vfi,\pifi,\lafi)$. Then we define $\vfi \ot \io_M \ot \io_N = \vfi \ot \io_{M \ot N}$ and $\lafi
\ot \io_M \ot \io_N = \lafi \ot \io_{M \ot N}$. Similarly for $\io_M \ot  \io_N \ot \vfi$.

In order to define $\io_M \ot \vfi \ot \io_N$, we use a permutation operator. Therefore let $\flip : M \ot L
\rightarrow L \ot M$ denote the flip $^*$-isomorphism, then we set $\io_M \ot \vfi \ot \io_N = (\vfi \ot \io_M \ot
\io_N) (\flip \ot \io_N)$. If $M$ acts on $H$ and $N$ acts on $K$, we define $\io_M \ot \lafi \ot \io_N = (\Sigma \ot
1_K)(\lafi \ot \io_M \ot \io_N)(\flip \ot \io_N)$ where $\Sigma : H_\vfi \ot H \rightarrow H \ot H_\vfi$ denotes the
flip transformation.

\bigskip\medskip

\section{The definition of a von Neumann algebraic quantum group}

\bigskip

In \cite{VaKust} we defined reduced \cst-algebraic quantum groups. But it is also possible to formulate the theory in
the von Neumann algebra setting as was done in \cite{VaKust2}. In this paper we will use this alternative approach
because it better suits our needs. It should be said however that there exists a natural bijection between these two
notions of quantum groups in the operator algebra setting.

\smallskip

\begin{definition} \label{prel2.def1}
Consider a von Neumann algebra and a unital normal $^*$-homomorphism $\de : M \rightarrow M \ot M$ such that $(\de \ot
\io)\de = (\io \ot \de)\de$. Assume moreover the existence of n.s.f. weights $\vfi$ and $\psi$ on $M$ such that
\begin{enumerate}
\item We have for all $a \in \cM_\vfi^+$ and $\om \in M^+_*$ that  $\vfi((\om \ot
\io)\de(a)) = \vfi(a) \, \om(1)$.
\item We have for all $a \in \cM_\psi^+$ and $\om \in M^+_*$ that  $\psi((\io \ot
\om)\de(a)) = \psi(a) \, \om(1)$.
\end{enumerate}
Then we call the pair $(M,\de)$ a von Neumann algebraic quantum group.
\end{definition}

\smallskip

An elaborate discussion about von Neumann algebraic quantum groups and their relation with their \cst-algebraic
counterparts can be found in \cite{VaKust2}.

\medskip

Let us  fix some further notations and terminology. Therefore consider a von Neumann algebraic quantum group $(M,\de)$.
Fix also a n.s.f. weight $\vfi$ on $M$ satisfying property (1) of definition \ref{prel2.def1} (such a weight is called
a left Haar weight of $(M,\de)$\,).

\smallskip

There exists a unique $\si$-strongly$^*$ closed mapping $S$ in $M$ such that
\begin{enumerate}
\item We have for all $a,b \in \cN_\vfi$ that $(\io \ot
\vfi)(\de(a^*)(1 \ot b)) \in D(S)$ and $S\bigl((\io \ot \vfi)(\de(a^*)(1 \ot b))\bigr) = (\io \ot \vfi)((1 \ot
a^*)\de(b))$
\item The space $\langle \, (\io \ot \vfi)(\de(a^*)(1 \ot b)) \mid a,b \in \cN_\vfi\,\rangle$ is a $\si$-strong$^*$ core for $S$.
\end{enumerate}

\smallskip

There exist a unique anti $^*$-automorphism $R$ on $M$ and a unique $\si$-strongly$^*$ continuous one parameter group
$\tau$ on $M$ such that $$R^2 = \io \hspace{2cm} R \tau_t = \tau_t R \ \ \ \text{ for all  }\  t \in \R \hspace{2cm} S
= R \tau_{-\frac{i}{2}} \ .$$ We call $S$ the antipode, $R$ the unitary antipode and $\tau$ the scaling group of our
quantum group $(M,\de)$. There exists a positive number $\nu > 0$ satisfying $\vfi\, \tau_t = \nu^{-t}\,\vfi$ for all
$t\in \R$. The number $\nu$ is referred to as the scaling constant of $(M,\de)$.

\medskip

Since $\flip(R \ot R) \de = \de R$, the n.s.f. weight $\psi = \vfi R$ satisfies condition (2) of definition
\ref{prel2.def1}.  We also get that $\psi\, \si_t^\vfi = \nu^{-t} \, \psi$ for $t \in \R$. Hence there exists a
positive self-adjoint operator $\sde$ affiliated to $M$ such that $\si^\vfi_t(\sde) = \nu^t\, \sde$ for all $t \in \R$
and $\psi
 = \vfi_\sde$ (for a precise definition of $\vfi_\sde$, see  definition 1.3 of \cite{Va2}). In other words, $\psi$
is absolutely continuous with respect to $\vfi$ and $\sde$ is the Radon Nykodim derivative of $\psi$ with respect to
$\vfi$. The element $\sde$ is called the modular element of $(M,\de)$.

\bigskip\medskip

\section{Fixing the data}
\label{data}

\bigskip

In this section we will introduce the data that we will use throughout this paper. So let $(M,\de_M)$ and $(N,\de_N)$
be two von Neumann algebraic quantum groups with left Haar weights $\vfi_M$ and $\vfi_N$ respectively. Without loss of
generality, we may and will  assume that $M$ and $N$ are in standard form with respect to Hilbert spaces $H_M$ and
$H_N$ respectively. We also fix GNS-constructions $(H_M,\io,\la_M)$ and $(H_N,\io,\la_N)$ for $\vfi_M$ and $\vfi_N$
respectively. Throughout this paper, the natural objects (like the antipode, the unitary antipode) associated to
$(M,\de_M)$ will get a sub- or superscript $M$. The same rule applies to $(N,\de_N)$.

\smallskip

Define the right Haar weights $\psi_M$ and $\psi_N$ on $(M,\de_M)$ and $(N,\de_N)$ respectively as $\psi_M =
\vfi_M\,R_M$ and $\psi_N = \vfi_N\,R_N$, where $R_M$ and $R_N$ are the unitary antipodes of  $(M,\de_M)$ and
$(N,\de_N)$  respectively.

\smallskip

Let $\sde_M$ and $\sde_N$ denote the modular elements of $(M,\de_M)$ and $(N,\de_N)$ respectively. So $\sde_M$  is the
Radon Nykodim derivative of $\psi_M$ with respect to $\vfi_M$ and $\sde_N$  is the Radon Nykodim derivative of $\psi_N$
with respect to $\vfi_N$. We define the GNS-constructions $(H_M,\io,\Gamma_M)$ and $(H_N,\io,\Gamma_N)$ of $\psi_M$ and
$\psi_N$ respectively by setting $\Gamma_M = (\la_M)_{\sde_M}$ and $\Gamma_N = (\la_N)_{\sde_N}$.

\bigskip

We also would like $(N,\de_N)$ to be some kind of\ \lq sub quantum group\rq\ of $(M,\de_M)$. Without trying to go
further into the somewhat delicate notion of a \ \lq sub quantum group\rq\,, we will impose a less restrictive
condition on the pair $M$,$N$.

\smallskip

So we will assume the existence of a normal injective $^*$-homomorphism $\al : M \rightarrow M \ot N$ such that
\begin{equation}
(\al \ot \io)\al = (\io \ot \de_N)\al \hspace{2cm} \text{and} \hspace{2cm} (\de_M \ot \io)\al = (\io \ot \al)\de_M \ .
\label{sub}
\end{equation}
The first equality means that $\al$ is a right action of $(N,\de_N)$ on $M$. These two conditions are equivalent to the
fact that there exists a morphism of quantum groups from $(M,\de_M)$ to $(N,\de_N)$, a fact that is proven in section
10 of \cite{Kust1}.

\medskip

We define the sub von Neumann algebra $Q$ of $M$ as $$Q = \{\,x \in M \mid \al(x) = x \ot 1 \,\} \ .$$ This $Q$ should
be thought of as a quantum analogue of the $L^\infty$-functions on the left coset space. The role of the
quasi-invariant measure on the left coset space will be played by any normal semi finite faithful weight on $Q$.
Therefore we fix some n.s.f. weight $\theta$ on $Q$ together with a GNS-construction $(H_\theta,\pi_\theta,\la_\theta)$
for it.

\smallskip

\begin{lemma} \label{data.lem1}
We have $\de_M(Q) \subseteq M \ot Q$.
\end{lemma}
\begin{proof}
Choose $x \in Q$. Take $y \in Q'$. Using equation (\ref{sub}), we get for all $\om \in M_*$, $$\al\bigl((\om \ot
\io)\de_M(x)\bigr) = (\om \ot \io \ot \io)\bigl((\io \ot \al)\de_M(x)\bigr) = (\om \ot \io \ot \io)\bigl((\de_M \ot
\io)\al(x)\bigr) = (\om \ot \io)(\de_M(x)) \ot 1 \ ,$$ implying that $(\om \ot \io)\de_M(x) \in Q$ and hence $(\om \ot
\io)(\de_M(x))\,y = y \,(\om \ot \io)(\de_M(x))$. We conclude that $\de_M(x)(1 \ot y) = (1 \ot y)\de_M(x)$. Therefore
$\de_M(x) \in (M' \ot Q')' = M \ot Q$.
\end{proof}

So we see, by the coassociativity of $\de_M$, that the mapping $\be : \pi_\theta(Q) \rightarrow M \ot \pi_\theta(Q)$,
defined such that $\be(\pi_\theta(x)) = (\io \ot \pi_\theta)\de_M(x)$ for all $x \in Q$, is a left action of
$(M,\de_M)$ on $Q$. In \cite{Va1}, Stefaan Vaes showed that any such a left action has an implementation by a unitary
corepresentation. So we define $\impl$ to be the unitary element in $M \ot B(H_\theta)$ such that $\impl^*$ is the
unitary implementation, as defined in definition 3.6 of \cite{Va1}, of the left action $\be$. Recall from \cite{Va1}
that $\impl$ enjoys the following properties:
\begin{eqnarray}
& & \hspace{-6.5cm} 1.\ \ (\io \ot \pi_\theta)\de_M(x) = \impl^*(1 \ot \pi_\theta(x))\impl \ \text{ for all } x \in Q \
, \label{impl1} \\ & & \hspace{-6.5cm} 2.\ \ (\de_M \ot \io)(\impl) = \impl_{13} \impl_{23} \ . \label{impl2}
\end{eqnarray}

\medskip\medskip

As a last piece of data we fix a unitary corepresentation $U$ of $(N,\de_N)$ on a Hilbert space $K$, i.e. $U$ is a
unitary element in $N \ot B(K)$ such that  $(\de_N \ot \io)(U) = U_{13}\,U_{23}$.

\bigskip\medskip

\section{The carrier Hilbert space of the induced corepresentation}
\label{carrier}

\bigskip

This section is devoted to the construction of the carrier Hilbert space of our induced corepresentation. The
construction of this Hilbert space is modeled on the classical case. We will also give an alternative description for
the tensor products of this carrier Hilbert space with any Hilbert space.

\smallskip

Throughout this section we fix a Hilbert space $H$.

\begin{definition} We define the subspace $\cP_H$ of $B(H) \ot M \ot B(K)$ as
$$\cP_H = \{\,X \in B(H) \ot M \ot B(K) \mid (\io \ot \al \ot \io)(X) = U_{34}^* X_{124} \,\} \ .$$
\end{definition}

It is easy to check that $\cP_H$ has the following multiplicative structure: we have for any $X \in \cP_H$ that
\begin{enumerate}
\item $(Y \ot 1)X \in \cP_H$ for all $Y \in B(H) \ot Q$\ ,
\item $X Y \in \cP_H$ for all $Y \in B(H) \ot Q \ot B(K)$ \ .
\end{enumerate}

\smallskip

Also the proof of the following result is very elementary.

\begin{lemma}
Let $X,Y \in \cP_H$, then $Y^* X$ belongs to $B(H) \ot Q \ot B(K)$.
\end{lemma}
\begin{proof}
By definition of $\cP_H$, we have that $$(\io \ot \al \ot \io)(Y^* X) = Y_{124}^* U_{34} U_{34}^* X_{124}^* = (Y^*
X)_{124} \ .$$ So we get for all $\om_1 \in B(H)_*$ and $\om_2 \in B(K)_*$ that $\al\bigl((\om_1 \ot \io \ot \om_2)(Y^*
X)\bigr) = (\om_1 \ot \io \ot \om_2)(Y^* X) \ot 1$ implying that $(\om_1 \ot \io \ot \om_2)(Y^* X) \in Q$. Arguing as
in the proof of lemma \ref{data.lem1}, the lemma follows.
\end{proof}

\medskip

Thanks to this lemma, we can now define a sesquilinear form $\langle \, ,  \rangle$ on $\cP_H \odot (H \ot H_\theta \ot
K)$ such that $$\langle X \ot v , Y \ot w \rangle = \langle (\io \ot \pi_\theta \ot \io)(Y^* X)\, v , w \rangle $$ for
all $X,Y \in \cP_H$ and  $v,w \in H \ot H_\theta \ot K$.

\begin{lemma}
The inproduct $\langle \, , \rangle$ on $\cP_H \odot (H \ot H_\theta \ot K)$ is positive.
\end{lemma}
\begin{proof}
Choose $X_1,\ldots,X_n \in \cP_H$ and $v_1,\ldots,v_n \in H \ot H_\theta \ot K$. Then $$\langle \, \sum_{i=1}^n X_i \ot
v_i , \sum_{i=1}^n X_i \ot v_i \, \rangle = \sum_{i,j=1}^n \langle (\io \ot \pi_\theta \ot \io)(X_j^* X_i) \, v_i , v_j
\rangle \ .$$ Define $T \in M_n(B(H) \ot Q \ot B(K))$ such that $T_{ji} = X_j^* X_i$ for all $i,j=1,\ldots\!,n$. It is
clear that $T \geq 0$. Now define $S \in M_n(B(H) \ot B(H_\theta) \ot B(K))$ by setting $S_{ji} = (\io \ot \pi_\theta
\ot \io)(T_{ji})$ for $i,j=1,\ldots\!,n$, then it is clear that $S \geq 0$. Thus, $$\langle \, \sum_{i=1}^n X_i \ot v_i
, \sum_{i=1}^n X_i \ot v_i \, \rangle = \sum_{i,j=1}^n \langle S_{ji}\,v_i , v_j \rangle \geq 0 \ .$$
\end{proof}

\medskip

Set $\cJ = \{\,x \in \cP_H \odot (H \ot H_\theta \ot K) \mid \langle x , x \rangle = 0 \,\}$. Then $(\cP_H \odot (H \ot
H_\theta \ot K))/\cJ$ carries the structure of a pre-Hilbert space in the usual way.

\begin{notation}
We define $\cK_H$ to be the completion of $(\cP_H \odot (H \ot H_\theta \ot K))/\cJ$, the inproduct on $\cK_H$ will
also be denoted by $\langle \, , \rangle$. For $X \in \cP_H$ and $v \in H \ot H_\theta \ot K$ we define the element $X
\otd v$ in $\cK_H$ as the equivalence class of $X \ot v$ in $(\cP_H \odot (H \ot H_\theta \ot K))/\cJ$.
\end{notation}

\smallskip

We only need to remember the following elementary facts about $\cK_H$:
\begin{enumerate}
\item The mapping $\cP_H \times (H \ot H_\theta \ot K) \rightarrow \cK_H : (X,v) \mapsto X \otd v$ is bilinear,
\item $\langle X \otd v , Y \otd w \rangle = \langle (\io \ot \pi_\theta \ot \io)(Y^* X) \, v , w \rangle$ for all $X,Y \in \cP_H$ and $v,w \in H \ot H_\theta \ot K$,
\item $\cK_H = [ \, X \otd v \mid X \in \cP_H, v \in H \ot H_\theta \ot K \,]$ \ .
\end{enumerate}

\medskip

\begin{remark} \rm
In the special case where $H = \C$, the Hilbert space $\cK_\C$ will turn out to be the carrier space of the induced
corepresentation. Therefore we set $\cP = \cP_\C$ and $\cK = \cK_\C$.
\end{remark}

\smallskip

It is also clear that the following two properties hold
\begin{enumerate}
\item Let $X \in \cP_H$ and $\om \in B(H)_*$, then $(\om \ot \io \ot \io)(X) \in \cP$.
\item Let $X \in \cP$, then $1_{B(H}) \ot X \in \cP_H$.
\end{enumerate}

\smallskip

In the next proposition we will establish a natural isomorphism between $\cK_H$ and $H \ot \cK$. But observe first the
following basic facts. Let $X$ be an element in $\cP_H$. Since the map $H \ot H_\theta \ot K \rightarrow \cK_H : v
\mapsto X \otd v$ is clearly bounded (with bound $\leq \|X\|$), the following holds:

Let $D \subseteq H \ot H_\theta \ot K$ such that $H \ot H_\theta \ot K = [D]$. Then$\cK_H = [\,X \otd v \mid X \in
\cP_H, v \in D  \,]$.

\begin{proposition}
There exists a unique unitary linear transformation $U_H : H \ot \cK \rightarrow \cK_H$ such that $U_H(v \ot (X \otd
w)) = (1 \ot X) \otd (v \ot w)$ for all $v \in H$, $X \in \cP$ and $w \in H_\theta \ot K$. Furthermore, the following
holds: Let $v \in H$, $(e_i)_{i \in I}$ an orthonormal basis of $H$, $X \in \cP_H$ and $w \in H_\theta \ot K$. Then
$\sum_{i \in I} \, \|(\om_{v,e_i} \ot \io \ot \io)(X) \otd w \|^2 < \infty$ and $$U_H^*(X \otd (v \ot w)) = \sum_{i \in
I} e_i \ot \bigl((\om_{v,e_i} \ot \io \ot \io)(X) \otd w\bigr) \ .$$
\end{proposition}
\begin{proof}
For $v_1,v_2 \in H$, $X_1,X_2 \in \cP$ and $w_1,w_2 \in H_\theta \ot K$, we have
\begin{eqnarray*}
& & \langle (1 \ot X_1) \otd (v_1 \ot w_1) , (1 \ot X_2) \otd (v_2 \ot w_2) \rangle
\\ & & \spat = \langle (\io \ot \pi_\theta \ot \io)((1 \ot X_2^*)(1 \ot X_1))\, (v_1 \ot w_1) , v_2 \ot w_2 \rangle
\\ & & \spat = \langle (1 \ot (\pi_\theta \ot \io)(X_2^* X_1))\, (v_1 \ot w_1) , v_2 \ot w_2 \rangle
\\ & & \spat =  \langle v_1 , v_2 \rangle \, \langle (\pi_\theta \ot \io)(X_2^* X_1) \, w_1 , w_2 \rangle
\\ &  & \spat = \langle v_1 , v_2 \rangle \, \langle X_1 \otd w_1 , X_2 \otd w_2 \rangle
= \langle v_1 \ot (X_1 \otd w_1) , v_2 \ot (X_2 \otd w_2) \rangle \ .
\end{eqnarray*}
This implies the existence of an isometric linear map $U_H : H \ot \cK \rightarrow \cK_H$ such that \newline $U_H(v \ot
(X \otd w)) = (1 \ot X) \otd (v \ot w)$ for all $v \in H$, $X \in \cP$ and $w \in H_\theta \ot K$.

\smallskip

Choose $X \in \cP_H$, $v \in H$ and $w \in H_\theta \ot K$. Then
\begin{eqnarray*}
& & \sum_{i \in I} \, \|(\om_{v,e_i} \ot \io \ot \io)(X) \otd w \|^2 = \sum_{i \in I}\,\langle (\pi_\theta \ot
\io)((\om_{v,e_i} \ot \io \ot \io)(X)^* (\om_{v,e_i} \ot \io \ot \io)(X))\, w , w \rangle  \\ & & \spat
\stackrel{(*)}{=} \langle (\pi_\theta \ot \io)((\om_{v,v} \ot \io \ot \io)(X^* X)) \, w, w \rangle = \langle (\io \ot
\pi_\theta \ot \io)(X^* X) \, (v \ot w) , v \ot w \rangle < \infty \ ,
\end{eqnarray*}
where we used the normality of $\io \ot \pi_\theta \ot \io$ and equation (\ref{prel.eq1}) in (*). Now we get for all $Y
\in \cP_H$, $v' \in H$ and $w' \in H_\theta \ot K$ that
\begin{eqnarray*}
& & \langle U_H\bigl(\,\sum_{i \in I} e_i \ot ((\om_{v,e_i} \ot \io \ot \io)(X) \otd w)\,\bigr) , Y \otd (v' \ot w')
\rangle
\\
& & \spat =   \sum_{i \in I}\,\langle U_H\bigl(\, e_i \ot ((\om_{v,e_i} \ot \io \ot \io)(X) \otd w)\,\bigr) , Y \otd (v
\ot w) \rangle \\ & & \spat =  \sum_{i \in I}\, \langle (1 \ot (\om_{v,e_i} \ot \io \ot \io)(X)) \otd (e_i \ot w) , Y
\otd (v' \ot w') \rangle
\\
& & \spat = \sum_{i \in I}\, \langle (\io \ot \pi_\theta \ot \io)\bigl(Y^*(1 \ot (\om_{v,e_i} \ot \io \ot
\io)(X))\bigr)\, (e_i \ot w) , v' \ot w' \rangle  \\ & & \spat = \sum_{i \in I}\, \langle (\pi_\theta \ot
\io)\bigl((\om_{v',e_i} \ot \io \ot \io)(Y)^*(1 \ot (\om_{v,e_i} \ot \io \ot \io)(X))\bigr)\, w ,  w' \rangle \\ & &
\spat = \langle (\pi_\theta \ot \io)((\om_{v,v'} \ot \io \ot \io)(Y^* X)) \, w , w' \rangle =  \langle (\io \ot
\pi_\theta \ot \io)(Y^* X) \, (v \ot w) , v' \ot w' \rangle \\ & & \spat = \langle X \otd (v \ot w) , Y \otd (v' \ot
w') \rangle \ .
\end{eqnarray*}
Hence $U_H\bigl(\,\sum_{i \in I} e_i \ot ((\om_{v,e_i} \ot \io \ot \io)(X) \otd w)\,\bigr) = X \otd (v \ot w)$.

\smallskip

From this all we can also conclude that $U_H$ has dense range implying that it is a unitary transformation.
\end{proof}

Let us now use this proposition to introduce the following notation:

\begin{notation}
For all  $X \in \cP_H$, we define the bounded linear operator $X_* : H \ot H_\theta \ot K \rightarrow H \ot \cK$ such
that $X_*\,v = U_H^*(X \otd v)$ for all $v \in  H \ot H_\theta \ot K$.
\end{notation}

So we get immediately the following basic facts
\begin{enumerate}
\item The mapping $\cP_H \rightarrow B(H \ot H_\theta \ot K, H \ot \cK) : X \rightarrow X_*$ is linear and isometric ,
\item $(Y_*)^*(X_*) = (\io \ot \pi_\theta \ot \io)(Y^* X)$ for all $X,Y \in \cP_H$ ,
\item $H \ot \cK = [ \, X_*\,v \mid  X \in \cP_H, v \in H \ot H_\theta \ot K \,]$.
\end{enumerate}

\smallskip

Notice that the definition of $U_H$ implies that $(1_{B(H)} \ot X)_* = 1_{B(H)} \ot X_*$ for all $X \in \cP$.

\smallskip

The previous proposition implies for $v \in H$, $(e_i)_{i \in I}$ an orthonormal basis of $H$, $X \in \cP_H$ and $w \in
H_\theta \ot K$ that $\sum_{i \in I} \|(\om_{v,e_i} \ot \io \ot \io)(X)_*\,w\| < \infty$ and
\begin{equation}
X_*\,(v \ot w) = \sum_{i \in I} e_i \ot (\om_{v,e_i} \ot \io \ot \io)(X)_*\, w \ . \label{car}
\end{equation}

\medskip

As to be expected, this map $\cP_H \rightarrow B(H \ot H_\theta \ot K, H \ot \cK) : X \rightarrow X_*$ also behaves
well with respect to the strong topology on bounded sets.

\begin{result} \label{car.res2}
Consider a bounded net $(X_i)_{i \in I}$ in $\cP_H$ and $X \in \cP_H$. Then
\begin{enumerate}
\item If $(X_i)_{i \in I} \rightarrow X$ strongly, then $((X_i)_*)_{i \in I} \rightarrow X_*$ strongly.
\item If $(X_i)_{i \in I} \rightarrow X$ strongly$^*$, then $((X_i)_*)_{i \in I} \rightarrow X_*$ strongly$^*$.
\end{enumerate}
\end{result}
\begin{proof}
\begin{trivlist}
\item[\ \,1)]  Since $(X_i)_{i \in I} \rightarrow X$ strongly,
the net $\bigl((X_i - X)^*(X_i-X)\bigr)_{i \in I}$ is a bounded net in $B(H) \ot Q \ot B(K)$ that converges to 0 in the
weak operator topology. Therefore the normality of $\io \ot \pi_\theta \ot \io$ implies that the net  $\bigl((\io \ot
\pi_\theta \ot \io)((X_i - X)^*(X_i-X))\bigr)_{i \in I}$ also converges to 0 in the weak operator topology. Thus
$\bigl(((X_i)_* - X_*)^*((X_i)_*-X_*))\bigr)_{i \in I}$ converges to 0 in the weak operator topology. It follows that
$((X_i)_*)_{i \in I}$ converges strongly to $X_*$.
\smallskip
\item[\ \,2)] Choose $Y \in \cP_H$. Since $(X_i^*)_{i \in I}  \rightarrow X^*$ strongly,  the net
$\bigl(\,(X_i - X)(X_i - X)^*\,\bigr)_{i \in I}$ is a bounded net that converges to 0 in the weak operator topology.
Therefore $\bigl(\,Y^* (X_i - X)(X_i - X)^* Y\,)_{i \in I}$ is a bounded net in  $B(H) \ot Q \ot B(K)$ that converges
to 0 in the weak operator topology. Thus the normality of $\io \ot \pi_\theta \ot \io$ implies that the net $\bigl((\io
\ot \pi_\theta \ot \io)(Y ^* (X_i - X)(X_i - X)^* Y)\bigr)_{i \in I}$ also converges to 0 in the weak operator
topology. In other words, the net $\bigl( [((X_i)_* - X_*)^* Y_*]^*\,[((X_i)_* - X_*)^* Y_*] \bigr)_{i \in I}$
converges to 0 in the weak operator topology, implying that $\bigl(((X_i)_*)^* Y_*\bigr)_{i \in I}$ converges strongly
to $(X_*)^* Y_*$. Hence $\bigl(((X_i)_*)^*( Y_* v)\bigr)_{i \in I}$ converges to $(X_*)^* (Y_* v)$ for all $v \in H \ot
H_\theta \ot K$. Because the net $\bigl(((X_i)_*)^*\bigr)_{i \in I}$ is bounded and $H \ot \cK = [ \, X_*\,v \mid  X
\in \cP_H, v \in H \ot H_\theta \ot K \,]$, we conclude from this all that $\bigl(((X_i)_*)^*\bigr)_{i \in I}$
converges strongly to $(X_*)^*$.
\end{trivlist}
\end{proof}

\medskip

We will also need the following elementary properties.

\begin{result} \label{car.res1}
Consider $X \in \cP_H$, then
\begin{enumerate}
\item $(X\,Y)_* = X_*\,(\io \ot \pi_\theta \ot \io)(Y)$ for all $Y \in B(H) \ot Q \ot B(K)$\ ,
\item $((a \ot 1 \ot 1)\,X)_* = (a \ot 1) \, X_*$ for all $a \in B(H)$ \ .
\end{enumerate}
\end{result}
\begin{proof}
\begin{trivlist}
\item[\ \,1.] Choose $v \in H \ot H_\theta \ot K$. We have for all $Z \in \cP_H$ and $w \in H \ot H_\theta \ot  K$ that
\begin{eqnarray*}
& & \langle (X\,Y)_* v , Z_* w \rangle = \langle (\io \ot \pi_\theta \ot \io)(Z^* X Y) \, v , w \rangle \\ & & \spat =
\langle (\io \ot \pi_\theta \ot \io)(Z^* X) (\io \ot \pi_\theta \ot \io)(Y) \, v , w \rangle = \langle X_* (\io \ot
\pi_\theta \ot \io)(Y) \, v , Z_*\,w \rangle \ ,
\end{eqnarray*}
from which it follows that $(X\,Y)_* v = X_* (\io \ot \pi_\theta \ot \io)(Y) \, v$.

\smallskip

\item[\ \,2.]  Choose $v_1,v_2 \in H$, $w_1,w_2 \in H_\theta \ot K$ and $Y \in \cP$, then
\begin{eqnarray*}
& & \langle \bigl((a \ot 1 \ot 1)X\bigr)_* \,(v_1 \ot w_1) , v_2 \ot Y_* \,w_1 \rangle = \langle (\io \ot \pi_\theta
\ot \io)((1 \ot Y)^* (a \ot 1 \ot 1)X) \, (v_1 \ot w_1) , v_2 \ot w_2 \rangle
\\ & & \spat = \langle (\io \ot \pi_\theta \ot \io)((1 \ot Y)^* X) \, (v_1 \ot w_1) , a^* v_2 \ot w_2 \rangle
= \langle X_* \,(v_1 \ot w_1) , a^* v_2 \ot Y_* \,w_2 \rangle
\\ & & \spat = \langle (a \ot 1)\,X_*\, (v_1 \ot w_1) , a^* v_2 \ot Y_*\,
w_2 \rangle \ ,
\end{eqnarray*}
implying that $((a \ot 1 \ot 1)\,X)_* = (a \ot 1) \, X_*$ for all $a \in B(H)$.
\end{trivlist}
\end{proof}

\bigskip\medskip

\section{The definition of the induced corepresentation}
\label{definition}

\bigskip

In this section we define the induced corepresentation as a partial isometry. In  a later section we prove the unitary
of this induced corepresentation under an extra (mild?) condition.

\medskip

Consider $X \in \cP$. Using equation (\ref{sub}), we see that
\begin{eqnarray*}
& & (\io \ot \al \ot \io)\bigl((\de_M \ot \io)(X)\bigr) = (\de_M \ot \io \ot \io)\bigl((\al \ot \io)(X)\bigr) \\ & &
\spat = (\de_M \ot \io \ot \io)(U_{23}^* X_{13}) = U_{34}^* (\de_M \ot \io)(X)_{124} \ .
\end{eqnarray*}
Consequently, $(\de_M \ot \io)(X)$ belongs to  $\cP_{H_M}$.

\medskip

\begin{proposition}
There exists a unique isometry $\lambda \in B(H_M \ot \cK)$  such that $$\lambda(v \ot X_* w)  = (\de_M \ot \io)(X)_*
\, \impl_{12}^* (v \ot w) $$ for all $v \in H_M$, $X \in \cP$ and $w \in H_\theta \ot K$. Moreover, $\lambda$ belongs
to $M \ot B(\cK)$.
\end{proposition}
\begin{proof}
Referring to equation (\ref{impl1})  we get for all $v_1,v_2 \in H_M$, $w_1,w_2 \in H_\theta \ot K$ and $X_1,X_2 \in
\cP$ that
\begin{eqnarray*}
& & \langle (\de_M \ot \io)(X_1)_*\,\impl_{12}^* (v_1 \ot w_1) , (\de_M \ot \io)(X_2)_*\,\impl_{12}^* (v_2 \ot w_2)
\rangle \\ & & \spat = \langle (\io \ot \pi_\theta \ot \io)((\de_M \ot \io)(X_2)^* (\de_M \ot \io)(X_1)) \impl_{12}^*
(v_1 \ot w_1) , \impl_{12}^* (v_2 \ot w_2) \rangle \\ & & \spat = \langle \impl_{12}\, (\io \ot \pi_\theta \ot
\io)((\de_M \ot \io)(X_2^* X_1))\, \impl_{12}^* (v_1 \ot w_1) ,
 v_2 \ot w_2 \rangle \\
& & \spat =  \langle (1 \ot (\pi_\theta \ot \io)(X_2^* X_1))\,(v_1 \ot w_1) ,
 v_2 \ot w_2 \rangle \\
& & \spat = \langle v_1,v_2 \rangle \, \langle (X_1)_* \, w_1 , (X_2)_* \, w_2 \rangle = \langle v_1 \ot (X_1)_*\,w_1 ,
v_2 \ot (X_2)_*\,w_2 \rangle \ .
\end{eqnarray*}
From this chain of equalities the existence of $\lambda$ follows in the usual way.
\medskip

Choose $a \in M'$. Take $v \in H_M$, $X \in \cP$ and $w \in H_\theta \ot K$. Then, applying result \ref{car.res1} twice
and remembering that $(\de_M \ot \io)(X) \in M \ot B(H_M) \ot B(K)$ and $\impl \in M \ot B(K)$, we get that
\begin{eqnarray*}
& & (\lambda(a \ot 1))(v \ot X_* w) =   \lambda ( a v \ot X^* w ) = (\de_M \ot \io)(X)_* \, \impl_{12}^* (av \ot w) \\
& & \spat  =  (\de_M \ot \io)(X)_* (a \ot 1 \ot 1) \, \impl_{12}^* (v \ot w) = [ (\de_M \ot \io)(X)(a \ot 1 \ot 1) ]_*
\, \impl_{12}^* (v \ot w) \\ & & \spat =  [ (a \ot 1 \ot 1)(\de_M \ot \io)(X) ]_* \, \impl_{12}^* (v \ot w) = (a \ot
1)\, (\de_M \ot \io)(X)_* \, \impl_{12}^* (v \ot w) \\ &  & \spat = ((a \ot 1)  \lambda) (v \ot X_* w) \ .
\end{eqnarray*}
Hence $\lambda(a \ot 1) = (a \ot 1)  \lambda$. We conclude from this that $\lambda$ belongs to $(M' \ot 1)'= M \ot
B(\cK)$.
\end{proof}

\medskip

The next proposition establishes that $\lambda^*$ is a corepresentation of $(M,\de_M)$ on $\cK$.

\begin{proposition}
We have that $(\de_M \ot \io)(\lambda) = \lambda_{23}\,\lambda_{13}$.
\end{proposition}
\begin{proof} Define $\de_M^{(2)} = (\de_M \ot \io)\de_M = (\io \ot \de_M)\de_M$. For all  $X \in \cP$, we have that
\begin{eqnarray*}
& &  (\io_{B(H \ot H)} \ot \al \ot \io)\bigl(\,(\de_M^{(2)} \ot \io)(X)\,\bigr) = (\io \ot \io \ot \al \ot \io)(\de_M
\ot \io \ot \io)(\de_M \ot \io)(X) \\ & & \spat = (\de_M \ot \io \ot \io \ot \io)(\io \ot \al \ot \io)(\de_M \ot
\io)(X) = (\de_M \ot \io \ot \io \ot \io)(U_{34}^* (\de_M \ot \io)(X)_{124}) \\ & & \spat = U_{34}^* \, (\de_M^{(2)}
\ot \io)(X)_{124} \ ,
\end{eqnarray*}
implying that $(\de_M^{(2)} \ot \io)(X)$ belongs to $\cP_{H_M \ot H_M}$.

\medskip

Let $W$ denote the multiplicative unitary of $(M,\de_M)$ in the GNS-construction $(H_M,\io,\la_M)$. We know that
$\de_M(x) = W^* (1 \ot x) W$ for all $x \in M$.

\smallskip

Choose $v_1,v_2  \in H_M$, $w \in H_\theta \ot K$ and $X \in \cP$. Fix also an orthonormal basis $(e_i)_{i \in I}$ for
$H_M$.
\begin{trivlist}
\item[\ \,1)] We have that
\begin{eqnarray}
(\de_M \ot \io)(\lambda)(v_1 \ot v_2 \ot X_* w) & = & (W_{12}^* \lambda_{23} W_{12})(v_1 \ot v_2 \ot X_* w) \nonumber
\\ & =  & (W_{12}^* \lambda_{23})(W(v_1 \ot v_2) \ot X_* w) \ . \label{corepeq1}
\end{eqnarray}

Choose $u_1,u_2 \in H_M$. Using equation (\ref{car}) in the first and last step of the next chain of equalities, we get
for all $p \in H_M$ and $q \in  H_\theta \ot K$ that
\begin{eqnarray*}
& & W_{12}^* (u_1 \ot (\de_M \ot \io)(X)_* (p \ot q)) \\ & & \ \  = \sum_{i \in I} W_{12}^*\bigl(u_1 \ot e_i \ot
(\om_{p,e_i} \ot \io \ot \io)\bigl((\de_M \ot \io)(X)\bigr)_* \,q\,\bigr) \\ & & \ \  = \sum_{i \in I} \sum_{j \in I}
W_{12}^*\bigl(e_j \ot e_i \ot [(\om_{u_1,e_j} \ot \om_{p,e_i} \ot \io \ot \io)\bigl(1 \ot (\de_M \ot \io)(X)\bigr)]_*
\, q \, \bigr) \\ & & \ \ = \sum_{i \in I} \sum_{j \in I} W_{12}^*\bigl(e_j \ot e_i \ot [(\om_{W^*(u_1\ot p),W^*(e_j
\ot e_i)}  \ot \io \ot \io)\bigl(W_{12}^*(1 \ot (\de_M \ot \io)(X))W_{12} \bigr)]_* \, q \, \bigr) \\ & & \ \ = \sum_{i
\in I} \sum_{j \in I} W^*(e_j \ot e_i) \ot  [(\om_{W^*(u_1\ot p),W^*(e_j \ot e_i)}  \ot \io \ot \io)\bigl((\de_M^{(2)}
\ot \io)(X) \bigr)]_* \, q \, \bigr) \\ & & \ \  = (\de_M^{(2)} \ot \io)(X)_* W_{12}^*(u_1 \ot p \ot q) \ .
\end{eqnarray*}
Hence $$W_{12}^* (u_1 \ot (\de_M \ot \io)(X)_* \,\impl_{12}^*(u_2 \ot w)) = (\de_M^{(2)} \ot \io)(X)_* W_{12}^*(u_1 \ot
\impl_{12}^*(u_2 \ot w)) \ .$$ This implies that
\begin{eqnarray*}
& & W_{12}^* \lambda_{23} (u_1 \ot u_2 \ot X_*w) = W_{12}^* (u_1 \ot (\de_M \ot \io)(X)_*\,\impl_{12}^*(u_2 \ot w)) \\
& & \spat = (\de_M^{(2)} \ot \io)(X)_* W_{12}^*(u_1 \ot \impl_{12}^*(u_2 \ot w)) = (\de_M^{(2)} \ot \io)(X)_* W_{12}^*
\impl_{23}^*(u_1 \ot u_2 \ot w)
\end{eqnarray*}

Combining this with equation (\ref{corepeq1}), we find that
\begin{eqnarray}
 (\de_M \ot \io)(\lambda)(v_1 \ot v_2 \ot X_* w) & = & (\de_M^{(2)} \ot \io)(X)_* W_{12}^* \impl_{23}^* W_{12} (v_1 \ot v_2 \ot w) \nonumber \\
&  = & (\de_M^{(2)} \ot \io)(X)_* (\de_M \ot \io)(\impl^*)_{123} (v_1 \ot v_2 \ot w) \nonumber \\ & = & (\de_M^{(2)}
\ot \io)(X)_* \,\impl_{23}^* \,\impl_{13}^* (v_1 \ot v_2 \ot w) \ ,        \label{corepeq2}
\end{eqnarray}
where we used equation (\ref{impl2}) in the last step.

\medskip

\item[\ \,2)] Let $\Sigma$ denote the flip map on $H_M \ot H_M$. Then we have that
\begin{eqnarray}
 (\lambda_{23} \,\lambda_{13})(v_1 \ot v_2 \ot X_* w)
& = & (\lambda_{23} \Sigma_{12} \lambda_{23} \Sigma_{12})(v_1 \ot v_2 \ot X_* w)  \nonumber \\ & = & (\lambda_{23}
\Sigma_{12} \lambda_{23} )(v_2 \ot v_1 \ot X_* w) \nonumber \\ & = & (\lambda_{23} \Sigma_{12})\bigl(v_2 \ot (\de_M \ot
\io)(X)_*\,\impl_{12}^*(v_1 \ot w)\bigr) \ . \label{corepeq3}
\end{eqnarray}

Choose $p \in H_M$, $q \in H_\theta \ot K$. Then equation (\ref{car}) implies that
\begin{eqnarray}
& & (\lambda_{23} \Sigma_{12})\bigl(v_2 \ot (\de_M \ot \io)(X)_*(p \ot q)\bigr) \nonumber \\ & & \spat = \sum_{i \in I}
(\lambda_{23} \Sigma_{12})(v_2  \ot e_i \ot (\om_{p,e_i} \ot \io \ot \io)\bigl((\de_M \ot \io)(X)\bigr)_* \,q\,)
\nonumber \\ & & \spat = \sum_{i \in I} \lambda_{23}(\,e_i \ot v_2 \ot (\om_{p,e_i} \ot \io \ot \io)\bigl((\de_M \ot
\io)(X)\bigr)_* \,q\,) \nonumber \\ & & \spat = \sum_{i \in I} e_i \ot (\de_M \ot \io)\bigl((\om_{p,e_i} \ot \io \ot
\io)\bigl((\de_M \ot \io)(X)\bigr)\bigr)_* \,\impl_{12}^*(v_2 \ot q)  \nonumber \\ & & \spat = \sum_{i \in I} e_i \ot
(\om_{p,e_i} \ot \io \ot \io \ot \io)\bigl((\de_M^{(2)} \ot \io)(X)\bigr)_* \, \impl_{12}^* (v_2 \ot q) \ .
\label{corepeq4}
\end{eqnarray}
Now,
\begin{trivlist}
\item[\ \,$\bullet$] For all $u \in H_M \ot H_\theta \ot K$,
\begin{eqnarray*}
& & \sum_{i \in I} \|(\om_{p,e_i} \ot \io \ot \io \ot \io)\bigl((\de_M^{(2)} \ot \io)(X)\bigr)_*  \,u\,\|^2 \\ & &
\spat = \sum_{i \in I} \langle (\io \ot \pi_\theta \ot \io)\bigl((\om_{p,e_i} \ot \io \ot \io \ot
\io)\bigl((\de_M^{(2)} \ot \io)(X)\bigr)^*
\\ & & \spat  \spat \spat (\om_{p,e_i} \ot \io \ot \io \ot \io)\bigl((\de_M^{(2)} \ot \io)(X)\bigr)\bigr) \, u ,u \rangle \\
& & \spat = \langle (\io \ot \pi_\theta \ot \io)\bigl( (\om_{p,p} \ot \io \ot \io \ot \io)\bigl((\de_M^{(2)} \ot
\io)(X^* X)\bigr)\bigr) \, u ,u \rangle \\ & & \spat \leq \|p\|^2 \, \|X\|^2 \, \|u\|^2 \ ,
\end{eqnarray*}
implying that the linear map $$H_M \ot H_\theta \ot K \rightarrow H_M \ot H_M \ot \cK : u \mapsto \sum_{i \in I} e_i
\ot (\om_{p,e_i} \ot \io \ot \io \ot \io)\bigl((\de_M^{(2)} \ot \io)(X)\bigr)_*\,  u$$ is bounded.
\smallskip
\item[\ \,$\bullet$] For $u_1 \in H_M$, $u_2 \ot H_\theta \ot K$, we get, by applying equation (\ref{car}) twice,
\begin{eqnarray*}
& & \sum_{i \in I} e_i \ot (\om_{p,e_i} \ot \io \ot \io \ot \io)\bigl((\de_M^{(2)} \ot \io)(X)\bigr)_*  (u_1 \ot u_2)
\\ & & \spat = \sum_{i \in I} \sum_{j \in I} e_i \ot e_j \ot (\om_{p,e_i} \ot \om_{u_1,e_j} \ot \io \ot
\io)\bigl((\de_M^{(2)} \ot \io)(X)\bigr)_* \, u_2 \\ & & \spat = (\de_M^{(2)} \ot \io)(X)_* (p \ot u_1 \ot u_2) \ .
\end{eqnarray*}
\end{trivlist}
Combining these fact, we get that $$\sum_{i \in I} e_i \ot (\om_{p,e_i} \ot \io \ot \io \ot \io)\bigl((\de_M^{(2)} \ot
\io)(X)\bigr)_* \, u = (\de_M^{(2)} \ot \io)(X)_* (p \ot u) $$ for all $u \in H_M \ot H_\theta \ot K$. Using this in
combination with the chain of equalities in (\ref{corepeq4}), we see that $$(\lambda_{23} \Sigma_{12})\bigl(v_2 \ot
(\de_M \ot \io)(X)_*(p \ot q)\bigr) = (\de_M^{(2)} \ot \io)(X)_* \,\impl_{23}^*(p \ot v_2 \ot q)\ .$$ Hence, by
equation (\ref{corepeq3}), $$(\lambda_{23} \lambda_{13})(v_1 \ot v_2 \ot X_* w)  = (\de_M^{(2)} \ot \io)(X)_*
\,\impl_{23}^* \,\impl_{13}^* (v_1 \ot v_2 \ot w) \ .$$
\end{trivlist}
Comparing the above equation with equation (\ref{corepeq2}), we  conclude that $$(\de_M \ot \io)(\lambda)(v_1 \ot v_2
\ot X_* w) = (\lambda_{23}\,\lambda_{13})(v_1 \ot v_2 \ot X_* w) \ .$$
\end{proof}

\medskip

As mentioned before, $\lambda$ is really the adjoint of the induced corepresentation itself:

\begin{notation}
We define $\rho = \lambda^*$. So $\rho$ is a surjective partial isometry in $M \ot B(\cK)$ such that \newline $(\de_M
\ot \io)(\rho) = \rho_{13} \, \rho_{23}$. We call $\rho$ the induced corepresentation associated to the quadruple
$(M,\de_M)$, \newline $(N,\de_N)$,$\al$,$U$ with respect to the GNS-construction $(H_\theta,\pi_\theta,\la_\theta)$. We
refer to $\cK$ as the carrier space of $\rho$.
\end{notation}

\bigskip\medskip

Our definition of the induced corepresentation (and its carrier space) depends on the choice of the n.s.f. weight on
$Q$, but it is no big surprise that all these different induced corepresentations are unitarily equivalent. Let us
quickly formalize this statement.

\smallskip

So let $\eta$ be another n.s.f. weight on $Q$ with GNS-construction $(H_\eta,\pi_\eta,\la_\eta)$. Then $u$ will denote
the canonical unitary transformation $u : H_\theta \rightarrow H_\eta$ (see e.g. equation (1) in section 3.16 of
\cite{Stramod}). Recall that $u \pi_\theta(x) u^* = \pi_\eta(x)$ for all $x \in Q$.

\smallskip

Let $\varrho$ denote the induced corepresentation of the quadruple $(M,\de_M)$, $(N,\de_N)$, $\al$, $U$ with respect to
the GNS-construction $(H_\eta,\pi_\eta,\la_\eta)$. Define $\cL$ to be the carrier space of $\varrho$

\smallskip

\begin{proposition}
There exists a unique unitary transformation $\cU : \cK \rightarrow \cL$ such that $\cU  X_* = X_* (u \ot 1)$ for all
$X \in \cP$. Let $H$ be any Hilbert space, then we have moreover that $(1 \ot \cU) X_*  = X_* (1 \ot u \ot 1)$ for all
$X \in \cP_H$.
\end{proposition}
\begin{proof}
For $X,Y \in \cP$ and $v,w \in H_\theta \ot K$, we have
\begin{eqnarray*}
& & \langle X_* (u \ot 1) v , Y_* (u \ot 1) w \rangle = \langle (\pi_\eta \ot \io)(Y^* X)  (u \ot 1) v ,   (u \ot 1) w
\rangle \\ & & \spat = \langle (u^* \ot 1)(\pi_\eta \ot \io)(Y^* X)  (u \ot 1) v , w \rangle =  \langle (\pi_\theta \ot
\io)(Y^* X)   v , w \rangle = \langle X_* v , Y_* w \rangle \ .
\end{eqnarray*}
This implies the existence of an isometry $\cU : \cK \rightarrow \cL$ such that $\cU(X_* v) = X_*(u \ot 1) v$ for all
$X \in \cP$ and $v \in H_\theta \ot K$. It is then immediately clear that $\cU$ has dense range and is therefore
unitary.

Let $H$ be any Hilbert space and $X$ an element in $\cP_H$. Choose $v \in H$, $w \in \ot H_\theta \ot K$. Take also an
orthonormal basis $(e_i)_{i \in I}$ for $H$, then
\begin{eqnarray*}
(1 \ot \cU)X_* (v \ot w) & = &  \sum_{i \in I} (1 \ot \cU)(e_i \ot (\om_{v,e_i} \ot \io \ot \io)(X)_* w) \\ & = &
\sum_{i \in I} e_i \ot (\om_{v,e_i} \ot \io \ot \io)(X)_* (u \ot 1) w \\ & = & X_* (v \ot (u \ot 1) w) = X_* (1 \ot u
\ot 1)(v \ot w) \ ,
\end{eqnarray*}
thus $(1 \ot \cU)X_* = X_* (1 \ot u \ot 1)$.
\end{proof}

\medskip

Let $\impla$ denote the adjoint of the unitary implementation of the left action $\gamma : \pi_\eta(Q) \rightarrow M
\ot \pi_\eta(Q)$, defined such that $\gamma(\pi_\eta(a)) = (\io \ot \pi_\eta)\de_M(x)$ for all $x \in Q$. By
proposition 4.1 of \cite{Va1}, we know that $\impla = (1 \ot u) \impl (1 \ot u^*)$. From this we easily infer that

\begin{proposition}
The corepresentations $\rho$ and $\varrho$ are unitarily equivalent, i.e. $\varrho = (1 \ot \cU) \rho (1 \ot \cU^*)$.
\end{proposition}
\begin{proof}
Choose $X  \in \cP$, $v \in H_M$, $w \in H_\theta \ot K$. Then the previous proposition implies that
\begin{eqnarray*}
& & ((1 \ot \cU) \rho)(v \ot X_* w) = (1 \ot \cU) (\de_M \ot \io)(X)_* \, \impl_{12}^* (v \ot w) \\ & & \spat = (\de_M
\ot \io)(X)_* (1 \ot u \ot 1)\impl_{12}^* (v \ot w) = (\de_M \ot \io)(X)_* \,\impla_{12}^* (1 \ot u \ot 1) (v \ot w) \\
& & \spat = \varrho (v \ot X_* (u \ot 1) w) = (\varrho (1 \ot \cU)) (v \ot X_* w) \ ,
\end{eqnarray*}
and the proposition follows.
\end{proof}

\bigskip\medskip

\section{The integrability condition and its consequences for the carrier space $\cK$}

\bigskip

Recall that the action $\al$ is called integrable if the set $\{\,x \in M^+ \mid \al(x) \in \cM_{\io \ot \vfi_N}^+\,\}$
is $\si$-weakly dense in $M^+$. In general, an action does not have to be integrable but we will show that our special
actions are integrable under a very mild condition (see proposition \ref{integr.prop3}).  We will show that if our
action $\al$ is integrable, it is possible to produce an extremely useful dense subset of the carrier space $\cK$.

\medskip

Let us first round up the usual suspects:
\begin{enumerate}
\item $\cM^+ := \{ \, x \in M^+ \mid \al(x) \in \cM_{\io \ot \vfi_N}^+ \, \}$, so $\cM^+$ is a hereditary cone in $M^+$,
\item $\cN := \{ \, x \in M \mid x^* x \in \cM^+ \, \}$, so $\cN$ is a left ideal in $M$ ,
\item $\cM := $ the linear span of $\cM^+$ , so $\cM$ is a subalgebra of $M$.
\end{enumerate}

\medskip\smallskip

\begin{lemma} \label{integr.lem3}
\begin{enumerate}
\item $(\om \ot \io)\de_M(x) \in \cM^+$  for all $x \in \cM^+$ and $\om \in M_*^+$,
\item $(\om \ot \io)\de_M(x) \in \cM$  for all $x \in \cM$ and $\om \in M_*$,
\item Let  $x \in \cN$ and $\om \in M_*$. Then $(\om \ot \io)\de_M(x) \in \cN$ and $$\|(\io \ot \la_N)\bigl(\al((\om \ot \io)\de_M(x))\bigr)\| \leq \|\om\|\,\|(\io \ot \la_N)(\al(x))\|$$  \end{enumerate}
\end{lemma}
\begin{proof}
\begin{trivlist}
\item [\ \,1)] Using equation (\ref{sub}) we get that
$$\al((\om \ot \io)\de_M(x)) = (\om \ot \io \ot \io)((\io \ot \al)\de_M(x)) = (\om \ot \io \ot \io)((\de_M \ot
\io)\al(x)) \ .$$ So we get for every $\eta \in M_*^+$ that $$(\eta \ot \io)\bigl(\al((\om \ot \io)\de_M(x))\bigr) =
(\, ( \om \ot \eta)\de_M \ot \io)(\al(x)) \ \in \cM_{\vfi_N}^+  \ .$$ It follows that $\al((\om \ot \io)\de_M(x)) \in
\cM_{\io \ot \vfi_N}^+$. Moreover, we have for all $\eta \in M_*^+$ that
\begin{eqnarray*}
\eta\bigl((\io \ot \vfi_N)\bigl(\al((\om \ot \io)\de_M(x))\bigr)\bigr) & = & \vfi_N\bigl((\eta \ot \io)\bigl(\al((\om
\ot \io)\de_M(x))\bigr)\bigr) \\ & =& (\om \ot \eta)\de_M((\io \ot \vfi_N)\al(x)) \ ,
\end{eqnarray*}
implying that $(\io \ot \vfi_N)\bigl(\al((\om \ot \io)\de_M(x))\bigr) = (\om \ot \io)\de_M( (\io \ot \vfi_N)\al(x)) $
and hence $$\|(\io \ot \vfi_N)\bigl(\al((\om \ot \io)\de_M(x))\bigr)\| \leq \|\om\| \, \|(\io \ot \vfi)\al(x)\| \ .$$
\smallskip
\item[\ \,2)] This follows immediately from 1.
\smallskip
\item[\ \,3)] Using the estimate $(\om \ot \io)(\de_M(x))^* ( \om \ot \io)(\de_M(x)) \leq \|\om\| \, \, (|\om| \ot \io)\de_M(x^*x) $, this is an easy consequence of 1.
\end{trivlist}
\end{proof}

\medskip

Now we prove that the integrability of $\al$ is easy to check.

\begin{proposition} \label{integr.prop3}
The action $\al$ is integrable $\Leftrightarrow$ There exists a non-zero element $x \in M^+$ such that $\al(x) \in
\cM_{\io \ot \vfi}^+$.
\end{proposition}
\begin{proof}
One implication is trivial. We will prove the other one. Therefore suppose that there exists a non-zero element $x \in
M^+$ such that $\al(x) \in \cM_{\io \ot \vfi}^+$. Let $\bar{\cN}$ denote the $\si$-weak closure of $\cN$ in $M$. Then
$\bar{\cN}$ is a $\si$-weakly closed left ideal in $M$ so there exists a projection $P$ in $M$ such that $\bar{\cN} = M
\, P$.

Choose $\om \in M_*$. By the previous lemma, we know that $(\om \ot \io)\de_M(y) \in \cN$ for every $y \in \cN$.
Therefore the normality of $(\om \ot \io)\de_M$ implies that $(\om \ot \io)\de_M(y) \in \bar{\cN}$ for all $y \in
\bar{\cN}$. In particular, we find that $(\om \ot \io)\de_M(P) \in \bar{\cN}$ which implies that $(\om \ot \io)\de_M(P)
=  (\om \ot \io)(\de_M(P)) \, P$.

From this all we conclude that $\de_M(P) (1 \ot P) = \de_M(P)$, thus $\de_M(P) \leq 1 \ot P$. Therefore lemma  6.4  of
\cite{VaKust} implies that $P = 0$ or $P=1$. But the assumption at the beginning of this proof tells us that $P \not=
0$, hence $P =1$ and $\bar{\cN} = M$.
\end{proof}

\medskip\medskip

For the rest of this paper we assume that our action $\al$ is integrable. Thus
\begin{enumerate}
\item $\cM^+$ is $\si$-weakly dense in $M^+$,
\item $\cM$ and $\cN$ are $\si$-weakly dense in $M$.
\end{enumerate}

\smallskip

Using Kaplansky's density theorem, this also implies that all these sets are strongly$^*$ dense in $M$ an this in a
bounded way, e.g. for every $x \in M$ there exists a net $(x_i)_{i \in I}$ in $\cM$ such that $\|x_i\|  \leq \|x\|$ for
all $i \in I$ and such that $(x_i)_{i \in I}$ converges strongly$^*$ to $x$.

\medskip

\begin{definition}
We define the linear map $T_\al : \cM \rightarrow M$ such that $T_\al(x) = (\io \ot \vfi_N)(\al(x))$ for all $x \in
\cM$.
\end{definition}

Remember that $T_\al(x) \geq 0$ for all $x \in \cM^+$. It is also clear that for all $x \in \cM$ and $a,b \in Q$, the
element $a x b$ belongs to $\cM$ and $T_\al(a x b) = a \, T_\al(x) \, b$.

\medskip

\begin{proposition}
The set $T_\al(\cM)$ is a $\si$-weakly dense two-sided $^*$-ideal of $Q$.
\end{proposition}
\begin{proof}
Let $x \in \cM$. Then $\al(x) \in \cM_{\io \ot \vfi_N}$. Thus $(\al \ot \io)\al(x) \in \cM_{\io \ot \io \ot \vfi_N}$
and $\al((\io \ot \vfi_N)\al(x)) = (\io \ot \io \ot \vfi_N)((\al \ot \io)\al(x)) $. By equation (\ref{sub}), we get
that $(\io \ot \de_N)\al(x) \in \cM_{\io \ot \io \ot \vfi_N}$ and $$(\io \ot \io \ot \vfi_N)((\io \ot \de_N)\al(x)) =
(\io \ot \io \ot \vfi_N)((\al \ot \io)\al(x))   = \al(T_\al(x)) \ .$$ By the left invariance of  $\vfi_N$, we know that
$(\io \ot \io \ot \vfi_N)((\io \ot \de_N)\al(x))  = (\io \ot \vfi_N)(\al(x)) \ot 1 = T_\al(x) \ot 1$, thus
$\al(T_\al(x)) = T_\al(x) \ot 1$, so $T_\al(x) \in Q$.

So we have proven that $T_\al$ is a two-sided ideal in $Q$. Using the techniques of the proof of proposition 2.5(1) in
\cite{Haa1}, we arrive at the conclusion that $T_\al(\cM)$ is $\si$-weakly dense  in $Q$.
\end{proof}

The techniques used in Chap.1, Sec.3, Cor.5 of Thm.2 of \cite{Dix1}, guarantee the following result (see also
proposition 5.2(2)  of \cite{Haa1}).

\begin{proposition} \label{integr.prop2}
There exists an increasing net $(x_i)_{i \in I}$ in $\cM^+$ such that $(T_\al(x_i))_{i \in I}$ converges strongly to 1.
\end{proposition}

\bigskip

So we see that the map $T_\al$ allows us to produce enough elements in $Q$. We will  generalize this construction to
produce enough elements in $\cP$. We borrowed the basic idea for this procedure from the classical theory of induced
group representations but have to use quite different techniques to obtain the relevant results. The starting point is
the following basic result.

\begin{proposition} \label{integr.prop1}
Consider $X \in M \ot B(K)$ such that $U_{23}(\al \ot \io)(X)$ belongs to $\cM_{\io \ot \vfi_N \ot \io}$. Then $(\io
\ot \vfi_N \ot \io)(U_{23} (\al \ot \io)(X))$ belongs to  $\cP$.
\end{proposition}
\begin{proof}
Because $U_{23} (\al \ot \io)(X)$ belongs to $\cM_{\io \ot \vfi_N \ot \io}$, we have that $(\al \ot \io \ot \io)(
U_{23} (\al \ot \io)(X))$ belongs to  $\cM_{\io \ot \io \ot \vfi_N \ot \io}$ and
\begin{equation}
(\io \ot \io \ot \vfi_N \ot \io)\bigl((\al \ot \io \ot \io)( U_{23} (\al \ot \io)(X))\bigr)  = (\al \ot \io)\bigl((\io
\ot \vfi_N \ot \io)(U_{23} (\al \ot \io)(X))\bigr) \ . \label{integr.eq1}
\end{equation}
On the other hand, since $U_{23} (\al \ot \io)(X) \in \cM_{\io \ot \vfi_N \ot \io}$, the left invariance of $\vfi_N$
implies that the element $(\io \ot \de_N \ot \io)(U_{23} (\al \ot \io)(X))$ belongs to $\cM_{\io \ot \io \ot \vfi_N \ot
\io}$ and $$(\io \ot \io \ot \vfi_N \ot \io)\bigl((\io \ot \de_N \ot \io)(U_{23} (\al \ot \io)(X))\bigr) = (\io \ot
\vfi_N \ot \io)(U_{23} (\al \ot \io)(X))_{13} \ .$$ Therefore $U_{24}^* (\io \ot \de_N \ot \io)(U_{23} (\al \ot
\io)(X))$ belongs to $\cM_{\io \ot \io \ot \vfi_N \ot \io}$ and
\begin{eqnarray}
& & (\io \ot \io \ot \vfi_N \ot \io)\bigl(U_{24}^* (\io \ot \de_N \ot \io)(U_{23} (\al \ot \io)(X))\bigr) \nonumber
\\ & & \spat = U_{23}^* \, (\io \ot \io \ot \vfi_N \ot \io)\bigl((\io \ot \de_N \ot \io)(U_{23} (\al \ot \io)(X))\bigr) \nonumber
\\ & & \spat = U_{23}^* \,(\io \ot \vfi_N \ot \io)(U_{23} (\al \ot \io)(X))_{13} \ . \label{integr.eq2}
\end{eqnarray}
Since $U$ is a unitary corepresentation, we get that
\begin{eqnarray*}
& &  U_{24}^* \, (\io \ot \de_N \ot \io)(U_{23}(\al \ot \io)(X))
\\ & & \spat =  U_{24}^* U_{24} U_{34} (\io \ot \de_N \ot \io)((\al \ot \io)(X))
=  U_{34} \,  (\io \ot \de_N \ot \io)((\al \ot \io)(X)) \ ,
\\ & & \spat =  U_{34}\,(\al \ot \io \ot \io)((\al \ot \io)(X)) = (\al \ot \io \ot \io)(U_{23}(\al \ot \io)(X)) \ .
\end{eqnarray*}
Combining this with equations (\ref{integr.eq1}) and (\ref{integr.eq2}), we see that $(\io \ot \vfi_N \ot \io)(U_{23}
(\al \ot \io)(X))$ belongs to $\cP$.
\end{proof}

\medskip

Now the natural question arises how to construct such elements $X$ mentioned in the previous proposition. This will be
dealt with in the next 3 results. First we introduce some terminology.

\smallskip

Define the norm continuous one-parameter group $\si^*$ on $N_*$ such that $\si^*_t(\om) = \om \, \si^{\psi_N}_t$ for
all $t \in \R$. Let $\om \in N_*$ and $z \in \C$. Remember that $\om \in D(\si^*_z)$ $\Leftrightarrow$ There exists
$\eta \in N_*$ such that $\om \si_z^{\psi_N} \subseteq \eta$. In the latter case, $\si^*_z(\om) = \eta$. Also note that
$\om \in D(\si^*_z)$ $\Leftrightarrow$ $\bar{\om} \in D(\si^*_{\bar{z}})$. If $\om \in D(\si^*_z)$, then
$\si^*_{\bar{z}}(\bar{\om}) = \overline{\si^*_z(\om)}$.

\smallskip

We denote the set of elements that are analytic with respect to $\si^*$ by $\cA$.

\medskip

\begin{lemma}
Let $\om \in D(\si^*_{\frac{i}{2}})$  and $\eta \in B(K)_*$. Then $(\io \ot \eta)([1 \ot (\om \ot \io)(U^*)]\,U^*)$
belongs to $D(\si^{\vfi_N}_{\frac{i}{2}})$ and $$\si^{\vfi_N}_{\frac{i}{2}}\bigl((\io \ot \eta)([1 \ot (\om \ot
\io)(U^*)]\,U^*)\bigr) = R\bigl((\io \ot \eta)([1 \ot (\,\si^*_{\frac{i}{2}}(\om)\, R \ot \io)(U)]\,U)\bigr) \ .$$
\end{lemma}
\begin{proof}
By proposition 6.8 of \cite{VaKust}, we get for every $t \in \R$ that
\begin{eqnarray*}
& & \si^{\vfi_N}_t\bigl((\io \ot \eta)([1 \ot (\om \ot \io)(U^*)]\,U^*)\bigr) = \si^{\vfi_N}_t\bigl((\io \ot \eta)((\io
\ot \om \ot \io)(U^*_{23} U^*_{13}))\bigr) \\ & & \spat = \si^{\vfi_N}_t\bigl((\io \ot \om \ot \eta)((\de_N \ot
\io)(U^*))\bigr) = \si^{\vfi_N}_t\bigl((\io \ot \om)\de_N((\io \ot \eta)(U^*))\bigr) \\ & & \spat = (\io \ot \om
\si^{\psi_N}_t)\bigl((\si^{\vfi_N}_t \ot \si^{\psi_N}_{-t})\de_N((\io \ot \eta)(U^*))\bigr) = (\io \ot
\si^*_t(\om))\bigl(\de_N(\tau^N_t((\io \ot \eta)(U^*)))\bigr) \ .
\end{eqnarray*}
Theorem 1.6 (4) of \cite{Wor5} tells us that $(\io \ot \eta)(U^*)$ belongs to $D(\tau^N_{\frac{i}{2}})$ and $$
\tau^N_{\frac{i}{2}}((\io \ot \eta)(U^*)) = R((\io \ot \eta)(U)) \ . $$ Combining these two facts, we see that $(\io
\ot \eta)([1 \ot (\om \ot \io)(U^*)]U^*)$ belongs to $D(\si^{\vfi_N}_{\frac{i}{2}})$ and
\begin{eqnarray*}
& & \si^{\vfi_N}_{\frac{i}{2}}\bigl((\io \ot \eta)([1 \ot (\om \ot \io)(U^*)]U^*)\bigr) = (\io \ot
\si^*_{\frac{i}{2}}(\om)\,)\bigl(\de_N(\tau^N_{\frac{i}{2}}((\io \ot \eta)(U^*)))\bigr)
\\ & & \spat = (\io \ot \si^*_{\frac{i}{2}}(\om)\,)\de_N\bigl(R((\io
\ot \eta)(U))\bigr) = (\io \ot \si^*_{\frac{i}{2}}(\om) \,)\bigl(\flip(R \ot R)\de_N((\io \ot \eta)(U))\bigr)
\\ & & \spat = R\bigl(\,(\,\si^*_{\frac{i}{2}}(\om)\, R \ot \io)\de_N((\io
\ot \eta)(U))\,\bigr) = R\bigl(\,(\io \ot \eta)\bigl((\,\si^*_{\frac{i}{2}}(\om)\, R \ot \io \ot \io)(U_{13}
U_{23})\bigr)\,\bigr)
\\ & & \spat = R\bigl((\io \ot \eta)([1 \ot (\,\si^*_{\frac{i}{2}}(\om)\, R \ot \io)(U)]\,U)\bigr) \ .
\end{eqnarray*}
\end{proof}

\medskip

Using the $^*$-operation we get the following variant of this result (needed for later purposes). Let $\om \in
D(\si^*_{-\frac{i}{2}})$ and $\eta \in B(K)_*$. Then $(\io \ot \eta)(U \,[1 \ot (\om \ot \io)(U)])$ belongs to
$D(\si^{\vfi_N}_{-\frac{i}{2}})$ and
\begin{equation} \si^{\vfi_N}_{-\frac{i}{2}}\bigl((\io \ot \eta)(U\,[1 \ot (\om \ot
\io)(U)])\bigr) = R\bigl((\io \ot \eta)(U^*\,[1 \ot (\,\si^*_{-\frac{i}{2}}(\om)\, R \ot \io)(U^*)])\bigr) \ .
\end{equation}

\medskip\smallskip

\begin{lemma} \label{integr.lem4}
Define the anti $^*$-automorphism $T : N \rightarrow N' : x \mapsto J_N x^* J_N$. Consider $\om \in
D(\si^*_{\frac{i}{2}})$ and $a \in \cN_{\vfi_N}$. Then the element $[a \ot (\om \ot \io)(U^*)]\,U^*$ belongs to
$\cN_{\vfi_N \ot \io}$ and $$(\la_N \ot \io)\bigl([a \ot (\om \ot \io)(U^*)]\,U^*\bigr) = (1 \ot
(\,\si^*_{\frac{i}{2}}(\om)\,R \ot \io)(U))\,(T R \ot \io)(U)\, (\la_N(a) \ot 1) \ .$$
\end{lemma}
\begin{proof}
Take an orthonormal basis $(e_i)_{i \in I}$ for $K$. Choose $v \in K$. Then the $\si$-weak lower semi-continuity of
$\vfi_N$ implies that
\begin{eqnarray*}
& & \vfi_N\bigl(\,(\io \ot \om_{v,v})\bigl(([a \ot (\om \ot \io)(U^*)]\,U^*)^*\, [a \ot (\om \ot \io)(U^*)]\,U^*\bigr)
\,\bigr)
\\ & & \spat = \vfi_N\bigl(\,\sum_{i \in I} (\io \ot \om_{v,e_i})([a \ot (\om \ot \io)(U^*)]\,U^*)^*\,
(\io \ot \om_{v,e_i})([a \ot (\om \ot \io)(U^*)]\,U^*)\,\bigr)
\\ & & \spat = \sum_{i \in I} \vfi_N\bigl(\,(\io \ot \om_{v,e_i})([a \ot (\om \ot \io)(U^*)]\,U^*)^*\, (\io \ot \om_{v,e_i})([a
\ot (\om \ot \io)(U^*)]\,U^*)\,\bigr)
\\ & & \spat = \sum_{i \in I} \|\la_N( (\io \ot \om_{v,e_i})([a \ot (\om \ot \io)(U^*)]\,U^*))\|^2
\\ & & \spat = \sum_{i \in I} \|\la_N( a\,(\io \ot \om_{v,e_i})([1 \ot (\om \ot \io)(U^*)]\,U^*))\|^2 \ ,
\end{eqnarray*}
so the previous lemma implies that
\begin{eqnarray*}
& & \vfi_N\bigl(\,(\io \ot \om_{v,v})\bigl(([a \ot (\om \ot \io)(U^*)]\,U^*)^*\,[a \ot (\om \ot \io)(U^*)]\,U^*\bigr)
\,\bigr)
\\ & & \spat = \sum_{i \in I} \|T\bigl(\si^{\vfi_N}_{\frac{i}{2}}(
(\io \ot \om_{v,e_i})([1 \ot (\om \ot \io)(U^*)]\,U^*))\bigr)\, \la_N(a)\|^2 \ ,
\\ & & \spat = \sum_{i \in I} \| (T R)\bigl((\io \ot \om_{v,e_i})([1 \ot (\,\si^*_{\frac{i}{2}}(\om)\,
R \ot \io)(U)]\,U)\bigr)\,\la_N(a) \|^2
\\ & & \spat = \sum_{i \in I} \|(\io \ot \om_{v,e_i})\bigl((1 \ot
(\,\si^*_{\frac{i}{2}}(\om)\, R \ot \io)(U))\, (T R \ot \io)(U)\bigl) \, \la_N(a) \|^2
\\ & & \spat = \sum_{i \in I} \langle\, (\io \ot \om_{v,e_i})\bigl((1 \ot (\,\si^*_{\frac{i}{2}}(\om
)\, R \ot \io)(U))\, (T R \ot \io)(U)\bigr)^*
\\ & & \spat \spat \spat \, (\io \ot \om_{v,e_i})\bigl((1 \ot (\,\si^*_{\frac{i}{2}}(\om)
\, R \ot \io)(U))\, (T R \ot \io)(U)\bigr) \la_N(a) , \la_N(a) \rangle
\\ & & \spat = \langle\, (\io \ot \om_{v,v})(\,[(1 \ot (\,\si^*_{\frac{i}{2}}(\om)\, R \ot \io)(U))\, (T R \ot \io)(U)]^*
\\ & & \spat \spat \ \ [ (1 \ot (\,\si^*_{\frac{i}{2}}(\om)\, R \ot \io)(U))\, (T R \ot \io)(U)] \, \la_N(a) , \la_N(a) \rangle
\ .
\end{eqnarray*}
By the lower semi-continuity of $\vfi_N$, this implies that
\begin{eqnarray*}
& & \vfi_N\bigl(\,(\io \ot \om)\bigl(([a \ot (\om \ot \io)(U^*)]\,U^*)^* [a \ot (\om \ot \io)(U^*)]\,U^*\bigr) \,\bigr)
\\ & & \spat = \langle\,(\io \ot \om)(\,[(1 \ot (\,\si^*_{\frac{i}{2}}(\om)\, R \ot \io)(U))\, (T R \ot \io)(U)]^* \\ & &
\spat \spat \ \  [ (1 \ot (\,\si^*_{\frac{i}{2}}(\om)\, R \ot \io)(U)) (T R \ot \io)(U)] \, \la_N(a) , \la_N(a) \rangle
\ .
\end{eqnarray*}
for all $\om \in B(K)_*$ from which we conclude that $[a \ot (\om \ot \io)(U^*)]\,U^*$ belongs to $\cN_{\vfi_N \ot
\io}$. By result \ref{prel.res1}, We have moreover for all $v \in K$ that
\begin{eqnarray*}
& & (\la_N \ot \io)([a \ot (\om \ot \io)(U^*)]\,U^*) \, v
\\ & & \spat = \sum_{i \in I} \la_N\bigl((\io \ot \om_{v,e_i})([a \ot (\om \ot \io)(U^*)]\,U^*)\bigr) \ot e_i
\\ & & \spat = \sum_{i \in I} \la_N\bigl(\,a\,(\io \ot \om_{v,e_i})([1 \ot (\om \ot \io)(U^*)]\,U^*)\bigr) \ot e_i
\\ & & \spat = \sum_{i \in I} T\bigl(\si^{\vfi_N}_{\frac{i}{2}}( (\io \ot \om_{v,e_i})([1 \ot (\om \ot \io)(U^*)]\,U^*))\bigr)\, \la_N(a)  \ot e_i
\\ & & \spat = \sum_{i \in I} (T R)\bigl((\io \ot \om_{v,e_i})([1 \ot (\,\si^*_{\frac{i}{2}}(\om)\, R \ot \io)(U)]\,U)\bigr)\,\la_N(a) \ot e_i
\\ & & \spat = \sum_{i \in I} (\io \ot \om_{v,e_i})\bigl((1 \ot (\,\si^*_{\frac{i}{2}}(\om)\, R \ot \io)(U))\, (T R
\ot \io)(U)\bigr)\, \la_N(a) \ot e_i
\\ & & \spat = (1 \ot (\,\si^*_{\frac{i}{2}}(\om)\,R
\ot \io)(U))\, (T R \ot \io)(U)\, (\la_N(a) \ot v) \ .
\end{eqnarray*}
\end{proof}

\medskip

This lemma implies easily the next one.

\begin{lemma} \label{integr.lem1}
Consider $\om \in D(\si^*_{\frac{i}{2}})$ and $X \in M \ot N$ such that $X \in \cN_{\io \ot \vfi_N}$. Then $[X  \ot
(\om \ot \io)(U^*)]\,U^*_{23}$ belongs to $\cN_{\io \ot \vfi_N \ot \io}$ and
\begin{eqnarray*}
& & (\io \ot \la_N \ot \io)([X \ot (\om \ot \io)(U^*)]\,U^*_{23})
\\ & & \spat = (1 \ot 1 \ot (\,\si^*_{\frac{i}{2}}(\om)\,  R \ot \io)(U))\,(TR \ot \io)(U)_{23}\, ((\io \ot \la_N)(X) \ot 1) \ .
\end{eqnarray*}
\end{lemma}
\begin{proof}
We know that there exists a bounded net $(X_i)_{i \in I}$ in $M \odot \cN_{\vfi_N}$ such that $(X_i)_{i \in I}$
converges strongly$^*$ to $X$ and $(\,(\io \ot \la_N)(X_i)\,)_{i \in I}$  is a bounded net that converges strongly$^*$
to $(\io \ot \la_N)(X)$. Then $\bigl(\,[X_i  \ot (\om \ot \io)(U^*)]\,U^*_{23}\,\bigr)_{i \in I}$ is surely a bounded
net that converges strongly$^*$ to $[X  \ot (\om \ot \io)(U^*)]\,U^*_{23}$. It follows easily from the previous lemma
that $[X_i  \ot (\om \ot \io)(U^*)]\,U^*_{23}$ belongs to $\cN_{\io \ot \vfi_N \ot \io}$ and
\begin{eqnarray*}
& & (\io \ot \la_N \ot \io)([X_i \ot (\om \ot \io)(U^*)]\,U^*_{23})
\\ & & \spat = (1 \ot 1 \ot (\,\si^*_{\frac{i}{2}}(\om)\,  R \ot \io)(U))\,(TR \ot \io)(U)_{23}\, ((\io \ot \la_N)(X_i) \ot 1) \ .
\end{eqnarray*}
Therefore the net $\bigl(\,(\io \ot \la_N \ot \io)([X_i \ot (\om \ot \io)(U^*)]\,U^*_{23})\,\bigr)_{i \in I}$ is
bounded and converges strongly$^*$ to $$(1 \ot 1 \ot (\,\si^*_{\frac{i}{2}}(\om)\,  R \ot \io)(U))\,(TR \ot
\io)(U)_{23}\, ((\io \ot \la_N)(X) \ot 1) \ .$$ Using the $\si$-strong$^*$-strong closedness of $\io \ot \la_N \ot
\io$, the lemma follows.
\end{proof}

\medskip\medskip

So we get easily the following two results that will be crucial to us.

\begin{result}
Consider $x \in \cN$ and $\om \in D(\si^*_{-\frac{i}{2}})$. Then $U_{23}(\al(x^*) \ot (\om \ot \io)(U)) \in \cN_{\io
\ot \vfi_N \ot \io}^*$.
\end{result}
\begin{proof}
The element $\bar{\om}$ belongs to $D(\si^*_{\frac{i}{2}})$. Hence the previous lemma implies that $[\al(x) \ot
(\bar{\om} \ot \io)(U^*)]\,U^*_{23}$ belongs to $\cN_{\io \ot \vfi_N \ot \io}$, implying that $U_{23}[\al(x^*) \ot (\om
\ot \io)(U)]$ belongs to $\cN_{\io \ot \vfi_N \ot \io}^*$.
\end{proof}

\medskip

\begin{remark} \rm
Now proposition \ref{integr.prop1} implies the following results (the second is a special case of the first one).
\begin{enumerate}
\item Consider $x \in \cN$ and $\om \in D(\si^*_{-\frac{i}{2}})$.
Let $X \in M \ot B(K)$ such that $(\al \ot \io)(X) \in \cN_{\io \ot \vfi_N \ot \io}$. Then $U_{23} (\al \ot \io)([x^*
\ot (\om \ot \io)(U)]\,X)$ belongs to $\cM_{\io \ot \vfi_N\ot \io}$ and the element $$(\io \ot \vfi_N \ot
\io)\bigl(U_{23} (\al \ot \io)([x^* \ot (\om \ot \io)(U)]\,X)\bigr)$$ belongs to $\cP$.
\item Consider $x \in \cM$ and $\om \in D(\si^*_{-\frac{i}{2}})$ and $y \in B(K)$ Then the element $U_{23}[\al(x) \ot (\io \ot \om)(U)\,y\,]$ belongs to $\cM_{\io \ot
\vfi_N\ot \io}$ and the element $$(\io \ot \vfi_N \ot \io)\bigl(U_{23}[\al(x) \ot (\om \ot \io)(U)\,y\,]\,\bigr)$$
belongs to $\cP$.
\end{enumerate}
\end{remark}

\medskip

\begin{result} \label{integr.res1}
Consider $x \in \cN$, $X \in \cP$ and $\om \in D(\si^*_{\frac{i}{2}})$. Then $(\al \ot \io)\bigl([x \ot (\om \ot
\io)(U^*)]\,X\bigr)$ belongs to $\cN_{\io \ot \vfi_N \ot \io}$.
\end{result}
\begin{proof}
By lemma \ref{integr.lem1} we get that $[\al(x) \ot (\om \ot \io)(U^*)]\, U_{23}^*$ belongs to $\cN_{\io \ot \vfi_N \ot
\io}$ implying that the element $(\al(x) \ot (\om \ot \io)(U^*)) U_{23}^* X_{13}$ belongs to $\cN_{\io \ot \vfi_N \ot
\io}$. Because $(\al \ot \io)(X) = U_{23}^* X_{13}$, the result follows.
\end{proof}

\medskip\smallskip

\begin{proposition}
There exists a directed set $I$ and nets $(a_i)_{i \in I}$ in $\cN$, $(\om_i)_{i \in I}$ in $\cA$ such that
\begin{enumerate}
\item $(\io \ot \vfi_N \ot \io)\bigl(U_{23} [\al(a_i^* a_i) \ot (\om_i \ot \io)(U) (\om_i \ot \io)(U)^*]\,
U_{23}^*\bigr) \leq 1$ for all $i \in I$.

\smallskip\smallskip

\item The net $\bigl(\,(\io \ot \vfi_N \ot \io)\bigl(U_{23} [\al(a_i^* a_i) \ot (\om_i \ot \io)(U) (\om_i \ot \io)(U)^*]\,
U_{23}^*\bigr)\,\bigr)_{i \in I}$ converges strongly to 1.
\end{enumerate}
\end{proposition}
\begin{proof}
Choose $b \in \cN$ and $\eta \in \cA$. Notice that $\bar{\eta} \in \cA$. Using lemma \ref{integr.lem1}, we see that
\begin{eqnarray}
& & (\io \ot \vfi_N \ot \io)\bigl(U_{23} [\al(b^* b) \ot (\eta \ot \io)(U)(\eta \ot \io)(U)^*] U_{23}^* \bigr)
\nonumber
\\ & & \spat = (\io \ot \la_N \ot \io)\bigl((\al(b) \ot 1)([1 \ot (\eta \ot \io)(U)^*] U^*)_{23}\bigr)^* \nonumber
\\ & & \spat \spat \ \  (\io \ot \la_N \ot \io)\bigl((\al(b) \ot 1)([1 \ot (\eta \ot \io)(U)^*] U^*)_{23}\bigr) \nonumber
\\ & & \spat = (\io \ot \la_N \ot \io)\bigl((\al(b) \ot 1)([1 \ot (\bar{\eta} \ot \io)(U^*)] U^*)_{23}\bigr)^* \nonumber
\\ & & \spat \spat \ \  (\io \ot \la_N \ot \io)\bigl((\al(b) \ot 1)([1 \ot (\bar{\eta} \ot \io)(U^*)] U^*)_{23}\bigr) \nonumber
\\ & & \spat = [(1 \ot 1 \ot (\,\si^*_{\frac{i}{2}}(\bar{\eta})\, R \ot \io)(U))\,(T R \ot \io)(U)_{23}\, ((\io \ot \la_N)(\al(b)) \ot 1)]^* \nonumber
\\ & & \spat\spat \ \ [(1 \ot 1 \ot (\,\si^*_{\frac{i}{2}}(\bar{\eta})\, R \ot \io)(U))\,(T R \ot \io)(U)_{23}\, ((\io \ot \la_N)(\al(b)) \ot 1)] \nonumber
\\ & & \spat = ((\io \ot \la_N)(\al(b))^* \ot 1)\, (T R \ot \io)(U^*)_{23} \nonumber
\\ & & \spat \spat \ \ (1 \ot 1 \ot (\,\si^*_{\frac{i}{2}}(\bar{\eta})\, R \ot \io)(U)^*\, (\,\si^*_{\frac{i}{2}}(\bar{\eta})\, R \ot \io)(U)\, )  \nonumber
\\ & & \spat \spat \ \ (T R \ot \io)(U)_{23} ((\io \ot \la_N)(\al(b)) \ot 1) \ , \label{integr.eq3}
\end{eqnarray}
Let us now get hold of some interesting elements:
\begin{trivlist}
\item[\ \,1.] Using proposition \ref{integr.prop2}, we get the existence of a net $(x_l)_{l \in L}$ in $\cM^+$
such that
\begin{itemize}
\item $(\io \ot \vfi_N)(\al(x_l)) \leq 1$ for all $l \in L$,
\item The net $\bigl(\,(\io \ot \vfi_N)(\al(x_l))\,\bigr)_{l \in L}$ converges strongly
to 1.
\end{itemize}
\medskip
\item[\ \,2.] Theorem 1.6 of \cite{Wor5} guarantees that the $^*$-algebra $\{\,(\phi \ot \io)(U)\mid \phi \in N^\sharp_*\,\}$ is a
non-degenerate sub$^*$-algebra of $B(K)$ ($N^\sharp_*$ is defined in definition 2.3 of \cite{VaKust2}). Hence
Kaplansky's density theorem implies the existence of a net $(\phi_j)_{j \in J}$ in $N^\sharp_*$ such that
\begin{itemize}
\item $\|(\phi_j \ot \io)(U)\| < 1$ for all $j \in J$,
\item $(\,(\phi_j \ot \io)(U)\,)_{j \in J}$ converges strongly$^*$ to 1.
\end{itemize}

Since $\{\, \si^*_{\frac{i}{2}}(\bar{\eta}) \, R  \, \mid \eta \in \cA \, \}$ is dense in $N_*$, this implies easily
the existence of a net $(\eta_p)_{p \in P}$ in $\cA$ such that

\begin{itemize}
\item $\| (\,\si^*_{\frac{i}{2}}(\bar{\eta}_p)\,R \ot \io)(U) \| \leq 1$ for all $p \in P$,
\item The net $\bigl(\,(\,\si^*_{\frac{i}{2}}(\bar{\eta}_p)\, R \ot \io)(U)\,\bigr)_{p \in P}$ converges strongly$^*$ to 1.
\end{itemize}
Thus we get that
\begin{itemize}
\item $(\,\si^*_{\frac{i}{2}}(\bar{\eta}_p)\,R \ot \io)(U)^* \, (\,\si^*_{\frac{i}{2}}(\bar{\eta}_p)\, R \ot \io)(U) \leq 1$ for all $p \in P$,
\item The net $\bigl(\,(\,\si^*_{\frac{i}{2}}(\bar{\eta}_p)\, R \ot \io)(U)^* \,
(\,\si^*_{\frac{i}{2}}(\bar{\eta}_p)\, R \ot \io)(U)\,\bigr)_{p \in P}$ converges strongly to 1.
\end{itemize}
\end{trivlist}

We will use these nets to construct the nets whose existence is claimed in the statement of the proposition. So let
$F(H_M \ot K)$ be the directed set of all finite subsets of $H_M \ot K$. Set $I = F(H_M \ot K) \times \N$ and put the
product ordering on this set.

\smallskip

Choose $i=(F,n) \in I$. We get first of all the existence of an element $l_i  \in L$ such that $$\|\, [(\io \ot
\vfi_N)(\al(x_{l_i}^*x_{l_i})) \ot 1)] \, v - v \| \leq \frac{1}{2n} $$ for all $v \in F$.

By the chain of equalities in (\ref{integr.eq3}), We have for $p \in P$ that
\begin{eqnarray*}
& & (\io \ot \vfi_N \ot \io)(U_{23}[\al(x_{l_i}^* x_{l_i}) \ot (\eta_p \ot \io)(U) (\eta_p \ot \io)(U)^*] U_{23}^*)
\\ & & \spat = ((\io \ot \la_N)(\al(x_{l_i}))^* \ot 1)\, (T R \ot \io)(U^*)_{23}
\\ & & \spat\ \ \ (1 \ot 1 \ot (\,\si^*_{\frac{i}{2}}(\bar{\eta}_p) \, R \ot \io)(U)^*\,(\, \si^*_{\frac{i}{2}}(\bar{\eta}_p)\, R \ot \io)(U)\, )
\\ & & \spat\  \ \ (T R \ot \io)(U)_{23} ((\io \ot \la_N)(\al(x_{l_i})) \ot 1)
\ .
\end{eqnarray*}
This implies that
\begin{itemize}
\item We have for $p \in P$ that
\begin{eqnarray*}
& & (\io \ot \vfi_N \ot \io)(U_{23}[\al(x_{l_i}^* x_{l_i}) \ot (\eta_p \ot \io)(U) (\eta_p \ot \io)(U)^*] U_{23}^*)
\\ & & \spat \leq ((\io \ot \la_N)(\al(x_{l_i}))^*
\ot 1) (T R  \ot \io)(U^*)_{23} (T R \ot \io)(U)_{23} ((\io \ot \la_N)(\al(x_{l_i})) \ot 1) \\ & & \spat = (\io \ot
\vfi_N)(\al(x_{l_i}^* x_{l_i})) \ot 1 \leq 1 \ .
\end{eqnarray*}
\item The net
$$\bigl(\,(\io \ot \vfi_N \ot \io)(U_{23}[\al(x_{l_i}^* x_{l_i}) \ot (\eta_p \ot \io)(U) (\eta_p \ot \io)(U)^*]
U_{23}^*) \,\bigr)_{p \in p}$$ converges strongly to $$ ((\io \ot \la_N)(\al(x_{l_i}))^* \ot 1) (T R  \ot
\io)(U^*)_{23} (T R \ot \io)(U)_{23} ((\io \ot \la_N)(\al(x_{l_i})) \ot 1) \ ,$$ which is equal to $(\io \ot
\vfi_N)(\al(x_{l_i}^* x_{l_i})) \ot 1$.
\end{itemize}

\smallskip

From this we infer the existence of an element $p_i \in P$ such that
\begin{eqnarray*}
& & \|\,(\io \ot \vfi_N \ot \io)(U_{23}[\al(x_{l_i}^* x_{l_i}) \ot (\eta_{p_i} \ot \io)(U) (\eta_{p_i} \ot \io)(U)^*]
U_{23}^*)
 \, v \\ & & \spat \spat \spat - ((\io \ot
\vfi_N)(\al(x_{l_i}^* x_{l_i})) \ot 1)\,v\| \leq \frac{1}{2n}
\end{eqnarray*}
for all $v \in F$.

Now set $a_i = x_{l_i}$ and $\om_i = \eta_{p_i}$. Then
\begin{enumerate}
\item $(\io \ot \vfi_N \ot \io)(U_{23} [\al(a_i^* a_i) \ot (\om_i \ot \io)(U) (\om_i \ot \io)(U)^*]
U_{23}^*) \leq 1$

\smallskip

\item For all $v \in F$,
$$\|\,(\io \ot \vfi_N \ot \io)(U_{23} [\al(a_i^* a_i) \ot (\om_i \ot \io)(U) (\om_i \ot \io)(U)^*] U_{23}^*)\,v  - v \|
\leq \frac{1}{n} \ . $$
\end{enumerate}
It is also clear that the last inequality implies that the net $$\bigl(\,(\io \ot \vfi_N \ot \io)(U_{23} [\al(a_i^*
a_i) \ot (\om_i \ot \io)(U) (\om_i \ot \io)(U)^*] U_{23}^*)\,\bigr)_{i \in I}$$ converges strongly to 1.
\end{proof}

\medskip

This proposition provides the necessary ammunition to prove the next lemma.

\begin{lemma} \label{integr.lem2}
The carrier space $\cK$ is the closed linear span of elements of the form $$(\io \ot \vfi_N \ot \io)\bigl(U_{23} (\al
\ot \io)([a^* \ot (\om \ot \io)(U)]\,X)\bigr)_* \, v $$ where $a \in \cN$, $\om \in \cA$ and $X \in M \ot B(K)$ such
that $(\al \ot \io)(X) \in \cN_{\io \ot \vfi_N \ot \io}$, $v \in H_\theta \ot K$.
\end{lemma}
\begin{proof}
Choose $Y \in \cP$ and $v \in H_\theta \ot K$. The previous lemma guarantees the exstence of a directed set $I$, a net
$(a_i)_{i \in I} \in \cN$ and a net $(\om_i)_{i \in I}$ in $\cA$ such that
\begin{itemize}
\item $(\io \ot \vfi_N \ot \io)(U_{23} [\al(a_i^* a_i) \ot (\om_i \ot \io)(U) (\om_i \ot
\io)(U)^*]\,U_{23}^*) \leq 1$ for all $i \in I$,
\item The net $$\bigl(\,(\io \ot \vfi_N \ot \io)(U_{23}[\al(a_i^* a_i) \ot (\om_i \ot
\io)(U) (\om_i \ot \io)(U)^*]\,U_{23}^*)\,\bigr)_{i \in I}$$ converges strongly to 1.
\end{itemize}
By result \ref{integr.res1}, we have for every $i \in I$ that $(\al \ot \io)([a_i \ot (\om_i \ot \io)(U)^*]Y)$ belongs
to $\cN_{\io \ot \vfi_N \ot \io}$ and
\begin{eqnarray*}
& & (\io \ot \vfi_N \ot \io)\bigl(U_{23} (\al \ot \io)((a_i^* \ot (\om_i \ot \io)(U))(a_i \ot (\om_i \ot \io)(U)^*)
Y)\bigr)
\\ & & \spat = (\io \ot \vfi_N \ot \io)(U_{23} [\al(a_i^* a_i) \ot (\om_i \ot
\io)(U)(\om_i \ot \io)(U)^*] U_{23}^* Y_{13})
\\ & & \spat = (\io \ot \vfi_N \ot \io)(U_{23} [\al(a_i^* a_i) \ot (\om_i \ot \io)(U)(\om_i \ot \io)(U)^*] U_{23}^*) \, Y \ .
\end{eqnarray*}
From this it follows that $$\bigl(\,(\io \ot \vfi_N \ot \io)\bigl(U_{23} (\al \ot \io)((a_i^* \ot (\om_i \ot
\io)(U))(a_i \ot (\om_i \ot \io)(U)^*) Y)\bigr)\,\bigr)_{i \in I}$$ is a bounded net that converges strongly to $Y$.
Referring to result \ref{car.res2}, we conclude from this that the net $$\bigl(\,(\io \ot \vfi_N \ot \io)\bigl(U_{23}
(\al \ot \io)((a_i^* \ot (\om_i \ot \io)(U))(a_i \ot (\om_i \ot \io)(U)^*) Y)\bigr)_*\,v\,\bigr)_{i \in I}$$ converges
to $Y_* v$.
\end{proof}

\medskip\smallskip

We will combine this result with the next elementary technical lemma on slice weights to get to the penultimate result
of this section. It will also be clear that the next lemma holds for any n.s.f. weight on any von Neumann algebra.

\begin{lemma}
Consider $X \in M \ot N \ot B(K)$ such that $X \in \cN_{\io \ot \vfi_N \ot \io}$. Then we have for every $\om \in
B(K)_*$ that $(\io \ot \io \ot \om)(X)$ belongs to $\cN_{\io \ot \vfi_N}$. Moreover, the following holds,

Let $v \in H_M$, $w \in K$ and $(e_i)_{i \in I}$ an orthonormal basis for $K$, then $$\sum_{i \in I} \|(\io \ot
\la_N)((\io \ot \io \ot \om_{w,e_i})(X))\, v \|  < \infty$$ and $$(\io \ot \la_N \ot \io)(X)\,(v \ot w) = \sum_{i \in
I} (\io \ot \la_N)((\io \ot \io \ot \om_{w,e_i})(X))\,v \ot e_i \ .$$
\end{lemma}
\begin{proof}
Let us first prove the first statement. Choose $Y \in \cM_{\io \ot \vfi_N \ot \io}^+$ and $\eta \in B(K)_*^+$. Then we
have for all $\phi \in M_*^+$ that $(\phi \ot \io)\bigl((\io \ot \io \ot \eta)(Y)\bigl) = (\phi \ot \io \ot \eta)(Y)$
which belongs to $\cM^+_{\vfi_N}$. This implies that $(\io \ot \io \ot \eta)(Y)$ belongs to $\cM_{\io \ot \vfi_N}^+$.

Let $\om \in B(K)_*$. Using the inequality $(\io \ot \io \ot \om)(X)^* (\io \ot \io \ot \om)(X) \leq \|\om\| \, (\io
\ot \io \ot |\om|)(X^* X) $, the above considerations imply that $(\io \ot \io \ot \om)(X)$ belongs to $\cN_{\io \ot
\vfi_N}$.

\smallskip

Now we turn to the second statement. Therefore choose $v \in H_M$, $w \in K$ and $(e_i)_{i \in I}$ an orthonormal basis
for $K$. We have that $$\sum_{i \in I} \|(\io \ot \la_N)((\io \ot \io \ot \om_{w,e_i})(X))\,v\|^2 = \sum_{i \in I}
\langle (\io \ot \vfi_N)\bigl((\io \ot \io \ot \om_{w,e_i})(X)^* (\io \ot \io \ot \om_{w,e_i})(X)\bigr)\, v , v \rangle
\ .$$ Since $\bigl(\,\sum_{i \in J} (\io \ot \io \ot \om_{w,e_i})(X)^* (\io \ot \io \ot \om_{w,e_i})(X)\,\bigr)_{J \in
F(I)}$ is an increasing net that converges strongly to $(\io \ot \io \ot \om_{w,w})(X^* X)$, the above equality and the
$\si$-weak lower semi-continuity of $\io \ot \vfi_N$ implies that $$ \sum_{i \in I} \|(\io \ot \la_N)((\io \ot \io \ot
\om_{w,e_i})(X))\,v\|^2 = \langle (\io \ot \vfi_N)((\io \ot \io \ot \om_{w,w})(X^* X))\,v , v \rangle  < \infty\ .$$

Take an orthonormal basis $(f_l)_{l \in L}$ of $H_M$. Then result \ref{prel.res1} implies that $$(\io \ot \la_N \ot
\io)(X)\,(v \ot w) = \sum_{(i,l) \in I \times L} f_l \ot \la_N((\om_{v,f_l} \ot \io \ot \om_{w,e_i})(X)) \ot e_i \ ,$$
On the other hand, referring to result \ref{prel.res1} once again, we get that
\begin{eqnarray*}
\sum_{i \in I} \, (\io \ot \la_N)((\io \ot \io \ot \om_{w,e_i})(X))\,v \ot e_i & = & \sum_{i \in I} \sum_{l \in L} f_l
\ot \la_N\bigl((\om_{v,f_l} \ot \io)((\io \ot \io \ot \om_{w,e_i})(X))\bigr) \ot e_i \\  \spat & = & \sum_{i \in I}
\sum_{l \in L} f_l \ot \la_N((\om_{v,f_l} \ot \io \ot \om_{w,e_i})(X)) \ot e_i \ .
\end{eqnarray*}
Comparing both expressions, we conclude that $$ (\io \ot \la_N \ot \io)(X)\,(v \ot w) = \sum_{i \in I} (\io \ot
\la_N)((\io \ot  \io \ot \om_{w,e_i})(X))\,v \ot e_i \ .$$
\end{proof}

\medskip

Now the proof of our last result is a mere formality.

\begin{proposition}
The carrier space $\cK$ is the closed linear span of elements of the form $$(\io \ot \vfi_N \ot \io)\bigl(U_{23}[\al(x)
\ot (\om \ot \io)(U)]\bigr)_*\,v \ ,$$ where $x \in \cM$, $\om \in \cA$  and $v \in H_\theta \ot K$.
\end{proposition}
\begin{proof}
Choose $a \in \cN$, $\om \in \cA$, $X \in M \ot B(K)$ such that $(\al \ot \io)(X) \in \cN_{\io \ot \vfi_N \ot \io}$, $u
\in H_\theta$ and $w\in K$. Lemma \ref{integr.lem2} tells us that $\cK$ is the closed linear span of elements of the
form $$(\io \ot \vfi_N \ot \io)\bigl(U_{23}(\al \ot \io)([a^* \ot (\om \ot \io)(U)]\,X)\bigr)_* (u \ot w) \ .$$

Take a basis $(e_i)_{i \in I}$ of $K$. There exists $\eta \in M_*^+$ such that $\eta(x) = \langle \pi_\theta(x) u , u
\rangle$ for all $x \in Q$. Since $M$ is in standard form, we can find an element $q \in H_M$ such that $\eta(x) =
\langle x q , q \rangle$ for all $x \in M$.

From the previous lemma, we know that $(\io \ot \om_{w,e_i})(X)$ belongs to $\cN$ for all $i \in I$ and that the net
$$\bigl(\, \sum_{i \in J} \, (\io \ot \la_N)\bigl(\al((\io \ot \om_{w,e_i})(X))\bigr)\,q \ot e_i\,\bigr)_{J \in F(I)}$$
converges to $(\io \ot \la_N \ot \io)((\al \ot \io)(X))(q \ot w)$. (*)

Fix $J \in F(I)$ for the moment. Then
\begin{eqnarray*}
& & \|\, (\io \ot \vfi_N \ot \io)(U_{23} (\al \ot \io)([a^* \ot (\om \ot \io)(U)]\,X))_* (u \ot w) \\ & & \spat \spat -
\sum_{i \in J} (\io \ot \vfi_N \ot \io)\bigl(U_{23}(\al(a^* (\io \ot \om_{w,e_i})(X)) \ot (\om \ot \io)(U))\bigr)_* (u
\ot e_i) \,\|^2
\\ & & \spat = \|\, (\io \ot \vfi_N \ot \io)(U_{23} (\al \ot \io)([a^* \ot (\om \ot \io)(U)]\,X)) (q \ot w)
\\ & & \spat \spat - \sum_{i \in J} (\io \ot \vfi_N \ot \io)\bigl(U_{23}(\al(a^* (\io \ot \om_{w,e_i})(X)) \ot (\om \ot \io)(U))\bigr) (q \ot e_i) \,\|^2
\\ & & \spat = \| (\io \ot \la_N \ot \io)([\al(a) \ot (\om \ot \io)(U)^*]\, U_{23}^*)^* (\io
\ot \la_N \ot \io)(X)(q \ot w)
\\ & & \spat \spat - \sum_{i \in J} (\io \ot \la_N \ot \io)([\al(a) \ot (\om \ot \io)(U)^*]
\,U_{23}^*)^* ((\io \ot \la_N)\bigl(\al((\io \ot \om_{w,e_i})(X))\bigr)\, q \ot e_i) \|^2 \ .
\end{eqnarray*}
Therefore the convergence in (*) implies that the net $$\bigl(\,\sum_{i \in J} (\io \ot \vfi_N \ot
\io)\bigl(U_{23}(\al(a^* (\io \ot \om_{w,e_i})(X)) \ot (\om \ot \io)(U))\bigr)_* (u \ot e_i)\,\bigr)_{J \in F(I)}$$
converges to $(\io \ot \vfi_N \ot \io)(U_{23}(\al \ot \io)([a^* \ot (\om \ot \io)(U)]\,X))_* (u \ot w)$. Therefore the
proposition follows from the considerations in  the beginning of the proof.
\end{proof}

\bigskip\medskip

\section{Unitarity of the induced corepresentation under the integrability condition}

\bigskip

Also in this section we will assume that $\al$ is integrable. The aim of this section is to prove the unitarity of the
induced corepresentation under this extra assumption.

\medskip

Define $\cN_0 = \langle \, (\om \ot \io)\de_M(a) \mid \om \in M_*, a \in \cN\,\rangle$. From lemma \ref{integr.lem3} we
know that $\cN_0 \subseteq \cN$ and that
\begin{equation} \label{unit.ineq1}
\|(\io \ot \la_N)\bigl(\al((\om \ot \io)\de_M(a))\bigr)\| \leq \|\om\|\, \|(\io \ot \la_N)(\al(a))\|
\end{equation}
for all $a \in \cN$ and $\om \in M_*$ (a fact that we will use several times in this section).

\medskip

Let us also define a certain closed subspace of $\cK$:

\smallskip

\begin{notation} We define $\cK_0$ as the closed linear span of the elements of the form
$$(\io \ot \vfi_N \ot \io)(U_{23} [\al(y^* x) \ot (\om \ot \io)(U)])_* \, v \ ,$$ where $x,y \in \cN_0$, $\om \in \cA$
and $v \in H_\theta\ot K$.
\end{notation}

\smallskip

We will ultimately prove that $\cK = \cK_0$.

\medskip\medskip

\begin{lemma} \label{unitarity.lem1}
We have that $\lambda(H_M \ot \cK) \subseteq H_M \ot \cK_0$.
\end{lemma}
\begin{proof}
Choose $x,y \in \cN$, $\om \in \cA$, $v \in H_M$ and $w \in H_\theta\ot K$. Then
\begin{eqnarray}
 & & \lambda\bigl(\,v \ot (\io \ot \vfi_N \ot \io)(U_{23}[\al(y^* x) \ot (\om \ot \io)(U)])_* \, w \,\bigr) \nonumber
\\  & & \spat = (\de_M \ot \io)\bigl((\io \ot \vfi_N \ot \io)(U_{23}[\al(y^* x) \ot (\om \ot \io)(U)])\bigr)_* \,\impl_{12}^*
 (v \ot w) \ . \label{unitarity.eq1}
\end{eqnarray}
Choose an orthonormal basis $(e_i)_{i \in I}$ for $H_M$. Choose $p \in H_M$ and $q \in H_\theta\ot K$. By equation
(\ref{car}), we know that the net
\begin{equation}
\bigl(\,\sum_{i \in J} e_i \ot (\om_{p,e_i} \ot \io \ot \io)\bigl((\de_M \ot \io )\bigl((\io \ot \vfi_N \ot
\io)(U_{23}[\al(y^* x) \ot (\om \ot \io)(U)])\bigr)\bigr)_* \, q\,\bigr)_{J \in F(I)} \label{unitarity.eq3}
\end{equation}
converges to $$(\de_M \ot \io)\bigl((\io \ot \vfi_N \ot \io)(U_{23}[\al(y^* x) \ot (\om \ot \io)(U)])\bigr)_* \,(p \ot
q) \ .$$

In the next part we will show that each of the sums in this net belongs to $H_M \ot \cK_0$. Therefore fix $j \in I$.
Because $U_{23}[\al(y^*) \ot (\om \ot \io)(U)]$ belongs to $\cN_{\io \ot \vfi_N \ot \io}^*$ and $\al(x) \ot 1$ belongs
to $\cN_{\io \ot \vfi_N \ot \io}$, we get that $U_{34}\,[(\de_M \ot \io )(\al(y^*)) \ot (\om \ot \io)(U)]$ belongs to
$\cN_{\io \ot \io \ot \vfi_N \ot \io}^*$ and that $(\de_M \ot \io)(\al(x)) \ot 1$ belongs to $\cN_{\io \ot \ot \io
\vfi_N \ot \io}$ and
\begin{eqnarray*}
& & (\de_M \ot \io)\bigl((\io \ot \vfi_N \ot \io)(U_{23} [\al(y^* x) \ot (\om \ot \io)(U)])\bigr)
\\ & & \spat = (\io \ot \io \ot \vfi_N \ot \io)\bigl(U_{34}\,[(\de_M \ot \io)(\al(y^*)) \ot (\om \ot \io)(U)]\,[(\de_M \ot \io)(\al(x)) \ot 1]\bigr)
\end{eqnarray*}
This implies that the net
\begin{eqnarray*}
& & \bigl(\,\,\sum_{l \in L} (\io \ot \vfi_N \ot \io)\bigl((\om_{e_l,e_j} \ot \io \ot \io \ot \io)(U_{34}\, [(\de_M \ot
\io)(\al(y^*)) \ot (\om \ot \io)(U)])
\\ & & \spat\spat\spat (\om_{p,e_l} \ot \io \ot \io \ot \io)((\de_M \ot
\io)(\al(x)) \ot 1)\,\bigr)\,\, \bigr)_{L \in F(I)}
\end{eqnarray*}
is a bounded net that converges to $$(\om_{p,e_j} \ot \io \ot \io)\bigl((\de_M \ot \io)\bigl((\io \ot \vfi_N \ot
\io)(U_{23}[\al(y^* x) \ot (\om \ot \io)(U)])\bigr)\bigr) \ .$$ Therefore result \ref{car.res2} guarantees that the net
\begin{eqnarray}
& & \bigl(\,\,\sum_{l \in L} (\io \ot \vfi_N \ot \io)\bigl((\om_{e_l,e_j} \ot \io \ot \io \ot \io)(U_{34}\, [(\de_M \ot
\io)(\al(y^*)) \ot (\om \ot \io)(U)]) \nonumber
\\ & & \spat\spat\spat (\om_{p,e_l} \ot \io \ot \io \ot \io)((\de_M \ot
\io)(\al(x)) \ot 1)\,\bigr)_*\,q\,\, \bigr)_{L \in F(I)} \label{unitarity.eq2}
\end{eqnarray}
converges to $$(\om_{p,e_j} \ot \io \ot \io)\bigl((\de_M \ot \io)\bigl((\io \ot \vfi_N \ot \io)(U_{23}[\al(y^* x) \ot
(\om \ot \io)(U)])\bigr)\bigr)_* \, q \ .$$ But, using equation (\ref{sub}), we have for all $l \in L$ that
\begin{eqnarray*}
& & (\io \ot \vfi_N \ot \io)\bigl(\,(\om_{e_l,e_j} \ot \io \ot \io \ot \io)(U_{34}\,[(\de_M \ot \io)(\al(y^*)) \ot (\om
\ot \io)(U)])  \\ & & \spat\spat\spat\spat  (\om_{p,e_l} \ot \io \ot \io \ot \io)((\de_M \ot \io)(\al(x)) \ot
1)\,\bigr)_*\,q
\\ & & \spat = (\io \ot \vfi_N \ot \io)\bigl(\,(\om_{e_l,e_j} \ot \io \ot \io \ot \io)(U_{34} \,[(\io \ot \al)(\de_M(y^*)) \ot (\om \ot \io)(U)])\,
\\ & & \spat \spat\spat\spat\spat\ \ \  (\om_{p,e_l} \ot \io \ot \io \ot \io)((\io \ot \al)(\de_M(x)) \ot 1)\,\bigr)_*\,q
\\ & & \spat = (\io \ot \vfi_N \ot \io)(U_{23} [\al\bigl((\om_{e_j,e_l} \ot \io)(\de_M(y))^* (\om_{p,e_l} \ot \io)(\de_M(x))\bigr) \ot (\om \ot \io)(U)]\,)_*\, q ,
\end{eqnarray*}
which clearly belongs to $\cK_0$. Therefore the convergence in (\ref{unitarity.eq2}) implies that $$(\om_{p,e_j} \ot
\io \ot \io)\bigl((\de_M \ot \io)\bigl((\io \ot \vfi_N \ot \io)(U_{23}[\al(y^* x) \ot (\om \ot \io)(U)])\bigr)\bigr)_*
\, q $$ belongs to $\cK_0$. Consequently, by referring to the convergence in (\ref{unitarity.eq3}), we see that the
element $$(\de_M \ot \io)\bigl((\io \ot \vfi_N \ot \io)(U_{23}[\al(y^* x) \ot (\om \ot \io)(U)])\bigr)_* (p \ot q) \
.$$

\smallskip

The lemma follows now from equation (\ref{unitarity.eq1}).
\end{proof}

\medskip

The next lemma is the last crucial step in the proof of the unitarity of the induced corepresentation.

\begin{lemma} \label{unit.lem1}
We have that $ [\,(a \ot 1) \lambda v \mid a \in M, v \in H_M \ot \cK_0\,] \subseteq H_M \ot \cK_0$.
\end{lemma}
\begin{proof}
Define the unitary element $V \in B(H_M) \ot M$ such that $V (\Ga_M(a) \ot \Ga_M(b)) = (\Ga_M \ot \Ga_M)(\de_M(a)(1 \ot
b))$ for all $a,b \in \cN_{\psi_M}$. Then the following holds
\begin{enumerate}
\item $(\om_{\Ga_M(a),\Ga_M(b)} \ot \io)(V^*) = (\psi_M \ot \io)(\de_M(b^*)(a \ot 1))$ for all $a,b \in \cN_{\psi_M}$,
\item $(\io \ot \de_M)(V) = V_{12} V_{13}$ ,
\item $\de_M(x) = V (x \ot 1)V^*$ for all $x \in M$.
\end{enumerate}
Call $A$ the norm closure of $\{\,(\om \ot \io)(V) \mid \om \in B(H_M)_*\,\}$, then $A$ is a $\si$-weakly dense
sub-\cst-algebra of $M$ such that $V$ belongs to the multiplier algebra $M(B_0(H_M) \ot A)$ (for once, the tensor
product is here the minimal  \cst-algebraic tensor product). These last properties follows from the fact that $V$ is a
manageable multiplicative unitary in the sense of \cite{Wor5}.

\medskip

First we prove some small technical results
\begin{trivlist}
\item[\ \,1.] Consider a Hilbert space $H$, $x \in \cN$, $v,w \in H_M$ and $a \in B(H_M \ot H)$. Then the element \newline $(\io \ot
\al)\bigl((\om_{v,x^* w} \ot \io \ot \io)(V_{13}^*(a \ot 1))\bigr)$ belongs to $\cN_{\io \ot \io \ot \vfi_N}$ and $$
\|(\io \ot \io \ot \la_N)\bigl((\io \ot \al)\bigl((\om_{v,x^* w} \ot \io \ot \io)(V_{13}^*(a \ot 1))\bigr)\bigr)\| \leq
\|v\|\,\|w\|\,\|a\|\,\|(\io \ot \la_N)(\al(x))\| \ .$$

\smallskip

If $u \in H_M$, we denote by $\theta_u$ the element in $B(\C,H_M)$ defined by $\theta_u(c) = c\,u$ for all $c \in \C$.
We have that
\begin{eqnarray*}
& & (\om_{v,x^* w} \ot \io \ot \io)(V_{13}^*(a \ot 1))^* \,(\om_{v,x^* w} \ot \io \ot \io)(V_{13}^*(a \ot 1))
\\ & & \spat = [(\theta_{x^* w}^* \ot 1 \ot 1) V_{13}^* (a \ot 1)(\theta_v \ot 1 \ot 1)]^* \, [(\theta_{x^* w}^* \ot 1 \ot 1) V_{13}^* (a
\ot 1)(\theta_v \ot 1 \ot 1)]
\\ & & \spat = (\theta_{v}^* \ot 1 \ot 1)(a^* \ot 1) V_{13} (\theta_{x^* w} \ot 1 \ot 1)(\theta_{x^* w}^*  \ot 1 \ot 1) V_{13}^* (a \ot 1)(\theta_v \ot 1 \ot 1)
\\ & & \spat = (\theta_{v}^* \ot 1 \ot 1)(a^* \ot 1) V_{13} (x^* \ot 1 \ot 1)
(\theta_w \theta_w^* \ot 1 \ot 1) (x \ot 1 \ot 1) V_{13}^* (a \ot 1)(\theta_v \ot 1\ot 1)
\\ & & \spat \leq \|w\|^2 \ (\theta_{v}^* \ot 1 \ot 1)(a^* \ot 1) V_{13} (x^* x \ot 1 \ot 1)  V_{13}^* (a \ot 1)(\theta_v \ot 1 \ot 1)
\\ & &  \spat =\|w\|^2 \  (\om_{v,v} \ot \io \ot \io)((a^* \ot 1)\de_M(x^* x)_{13}(a \ot 1)) \ .
\end{eqnarray*}
Therefore
\begin{eqnarray*}
& & (\io \ot \al)\bigl((\om_{v,x^* w} \ot \io \ot \io)(V_{13}^*(a \ot 1))\bigr)^*\, (\io \ot \al)\bigl((\om_{v,x^* w}
\ot \io \ot \io)(V_{13}^*(a \ot 1))\bigr)
\\ & & \spat \leq \|w\|^2 \  (\om_{v,v} \ot \io \ot  \io \ot \io)((a^* \ot 1 \ot 1)(\io \ot \al)(\de_M(x^* x))_{134}(a \ot 1 \ot 1))
\\ & & \spat = \|w\|^2 \ (\om_{v,v} \ot \io \ot  \io \ot \io)((a^* \ot 1 \ot 1)(\de_M \ot \io)(\al(x^* x))_{134}(a \ot 1 \ot 1))  \ .
\end{eqnarray*}
Since $\al(x) \in \cN_{\io \ot \vfi_N}$, this implies immediately that the element $$(\io \ot \al)\bigl((\om_{v,x^* w}
\ot \io \ot \io)(V_{13}^*(a \ot 1))\bigr)^* \, (\io \ot \al)\bigl((\om_{v,x^* w} \ot \io \ot \io)(V_{13}^*(a \ot
1))\bigr)$$ belongs to $\cM_{\io \ot \io \ot \vfi_N}^+$ and
\begin{eqnarray*}
& & \hspace{-0.5cm} (\io \ot \io \ot \vfi_N)\bigl(\,(\io \ot \al)\bigl((\om_{v,x^* w} \ot \io \ot \io)(V_{13}^*(a \ot
1))\bigr)^*\,(\io \ot \al)\bigl((\om_{v,x^* w} \ot \io \ot \io)(V_{13}^*(a \ot 1))\bigr)\,\bigr)
\\ & & \hspace{-0.5cm} \spat \leq \|w\|^2 \ (\io \ot \io \ot \vfi_N)\bigl(\,(\om_{v,v} \ot \io \ot \io \ot \io)((a^* \ot 1 \ot 1)(\de_M \ot \io)(\al(x^* x))_{134}(a \ot 1 \ot 1))\,\bigr) \\ & & \hspace{-0.5cm} \spat = \|w\|^2 \ (\om_{v,v} \ot \io \ot \io)((a^* \ot 1)\de_M((\io \ot
\vfi_N)\al(x^* x))_{13} (a \ot 1)) \ ,
\end{eqnarray*}
implying that the element $(\io \ot \al)\bigl((\om_{v,x^* w} \ot \io \ot \io)(V_{13}^*(a \ot 1))\bigr)$ belongs to
$\cN_{\io \ot \io \ot \vfi_N}$ and
\begin{eqnarray*}
& &  \|(\io \ot \io \ot \la_N)\bigl((\io \ot \al)\bigl((\om_{v,x^* w} \ot \io \ot \io)(V_{13}^*(a \ot
1))\bigr)\bigr)\|^2  \\ & & \spat \leq \|v\|^2 \,\|w\|^2 \,\|a\|^2 \,\|(\io \ot \vfi_N)(\al(x^* x))\| = \|v\|^2
\,\|w\|^2 \,\|a\|^2 \,\|(\io \ot \la_N)(\al(x)))\|^2 \ .
\end{eqnarray*}

\medskip

\item[\ \,2.] Consider $x \in \cN$ and $v,w \in H_M$. Then $(\om_{v,x^* w} \ot \io)(V^*)$ belongs to $\cN$ and
$$\|(\io \ot \la_N)\bigl(\al((\om_{v,x^* w} \ot \io)(V^*))\bigr)\| \leq \|v\|\,\|w\|\,\|(\io \ot \la_N)(\al(x))\| \ .$$

\smallskip

This is just a special case of the result proven in the previous part.

\medskip

\item[\ \,3.] Let $x \in \cN$. Then $(\de_M \ot \io)\al(x) \in \cN_{\io \ot \io \ot \vfi_N}$ and
$$\|(\io \ot \io \ot \la_N)((\de_M \ot \io)\al(x))\|  = \|(\io \ot \la_N)(\al(x))\| \ .$$

\smallskip

We have that $(\de_M \ot \io)(\al(x))^* (\de_M \ot \io)(\al(x)) = (\de_M \ot \io)\al(x^* x)$. Because $\al(x^* x)$
belongs to $\cM_{\io \ot \vfi_N}^+$, this implies that $(\de_M \ot \io)(\al(x))^* (\de_M \ot \io)(\al(x))$ belongs to
$\cM_{\io \ot \io \ot \vfi_N}^+$ and $$\|(\io \ot \io \ot \vfi_N)\bigl(\,(\de_M \ot \io)(\al(x))^* (\de_M \ot
\io)(\al(x))\,\bigr)\| = \|\de_M((\io \ot \vfi_N)\al(x^* x))\| = \|(\io \ot \vfi_N)\al(x^* x)\| \ .$$

\end{trivlist}

Having dealt with these elementary technical issues, we can start with the essential part of the proof. Let $\cT$
denote the Tomita-algebra of $\psi_N$.  By properties 2 and 3 above, we get that
\begin{eqnarray*}
& & [\,(\io \ot \io \ot \la_N)\bigl((\de_M \ot \io)\al((\om_{v,x^* w} \ot \io)(V^*))\bigr) \mid x \in \cN,v,w \in H_M
\,]
\\ & & \spat = [\,(\io \ot \io \ot \la_N)\bigl((\de_M \ot \io)\al((\om_{\Ga_M(b c),x^* \Ga_M(a)} \ot
\io)(V^*))\bigr) \mid x \in \cN, a,b \in \cN_{\psi_M}, c \in \cT \,]
\\ & & \spat = [\,(\io \ot \io \ot \la_N)\bigl((\de_M \ot \io)\al((\psi_M \ot \io)(\de_M(a^* x)(b c \ot 1)))\bigr)
 \mid x \in \cN, a,b \in \cN_{\psi_M}, c \in \cT \,]
\\ & & \spat = [\,(\io \ot \io \ot \la_N)\bigl((\de_M \ot \io)\al((\psi_M \ot \io)((\si^{\psi_N}_i(c)
\ot 1) \de_M(a^* x)(b  \ot 1)))\bigr) \mid
\\ & & \spat\spat  x \in \cN, a,b \in \cN_{\psi_M}, c \in \cT \,]
\\ & & \spat \subseteq [\,(\io \ot \io \ot \la_N)((\de_M \ot \io)\al(y)) \mid y \in \cN_0 \, ] \ .
\end{eqnarray*}
Therefore,
\begin{eqnarray*}
& & [\,(\io \ot \io \ot \la_N)((\de_M \ot \io)\al(y))\,(a \ot 1)  \mid y \in \cN_0, a \in A \, ]
\\ & & \spat \supseteq [\,(\io \ot \io \ot \la_N)\bigl((\de_M \ot \io)\al((\om_{v,x^* w} \ot \io) (V^*))\bigr)\,(a \ot 1) \mid x \in \cN,v,w \in H_M, a \in A\,]
\\ & & \spat \supseteq [\,(\io \ot \io \ot \la_N)\bigl((\de_M \ot \io)\al((\om_{y v,x^* w} \ot \io) (V^*))\bigr)\,(a \ot 1) \mid
\\ & & \spat\spat x \in \cN, y \in B_0(H_M), v,w \in H_M, a \in A\,]
\\ & & \spat \supseteq [\,(\io \ot \io \ot \la_N)\bigl((\io \ot \al)\de_M((\om_{y v,x^* w} \ot \io) (V^*))\bigr)\,(a \ot 1) \mid
\\ & & \spat\spat x \in \cN, y \in B_0(H_M), v,w \in H_M, a \in A\,]
\\ & & \spat = [\,(\io \ot \io \ot \la_N)\bigl((\io \ot \al)((\om_{y v,x^* w} \ot
\io \ot \io)(V_{13}^* V_{12}^*))\bigr)\,(a \ot 1) \mid
\\ & & \spat \spat  x \in \cN, y \in B_0(H_M), v,w \in H_M, a \in A\,]
\\ & & \spat = [\,(\io \ot \io \ot \la_N)\bigl((\io \ot \al)((\om_{v,x^* w} \ot
\io \ot \io)(V_{13}^* [V^*(y \ot a)]_{12}))\bigr) \mid
\\ & & \spat\spat x \in \cN, y \in B_0(H_M), v,w \in H_M, a \in A\,] \ .
\end{eqnarray*}
Referring to property 1 above, we infer from this that
\begin{eqnarray*}
& & [\,(\io \ot \io \ot \la_N)((\de_M \ot \io)\al(y))\,(a \ot 1)  \mid y \in \cN_0, a \in A \, ] \\ & & \spat
\supseteq [\,(\io \ot \io \ot \la_N)\bigl((\io \ot \al)((\om_{v,x^* w} \ot \io \ot \io)(V_{13}^* \,(V^* z)_{12}))\bigr)
\mid
\\ & & \spat\spat x \in \cN, z \in B_0(H_M) \ot A, v,w \in H_M\,]
\\ & & \spat = [\,(\io \ot \io \ot \la_N)\bigl((\io \ot \al)((\om_{v,x^* w} \ot
\io \ot \io)(V_{13}^* \,z_{12}))\bigr) \mid
\\ & & \spat\spat x \in \cN, z \in B_0(H_M) \ot A, v,w \in H_M\,]
\\ & & \spat \supseteq [\,(\io \ot \io \ot \la_N)\bigl((\io \ot \al)((\om_{v,x^* w} \ot
\io \ot \io)(V_{13}^* \,(y \ot a)_{12}))\bigr) \mid
\\ & & \spat\spat x \in \cN, y \in B_0(H_M), v,w \in H_M, a \in A\,]
\\ & & \spat = [\,(\io \ot \io \ot \la_N)\bigl((\io \ot \al)(a \ot (\om_{y v,x^* w} \ot
\io)(V^*) )\bigr) \mid
\\ & & \spat\spat x \in \cN, y \in B_0(H_M), v,w \in H_M, a \in A\,]
\\ & &  \spat  = [\, a \ot  (\io \ot \la_N)\bigl(\al((\om_{y v,x^* w} \ot
\io)(V^*) )\bigr) \mid x \in \cN, y \in B_0(H_M), v,w \in H_M, a \in A\,] \ .
\end{eqnarray*}
Hence, invoking property 2 above, we conclude that
\begin{eqnarray*}
& & [\,(\io \ot \io \ot \la_N)((\de_M \ot \io)\al(y))\,(a \ot 1)  \mid y \in \cN_0, a \in A \, ]
\\  & & \spat \supseteq [\, a \ot  (\io \ot \la_N)\bigl(\al((\om_{v,x^* w} \ot
\io)(V^*) )\bigr) \mid x \in \cN, v,w \in H_M, a \in A\,]
\\ & &  \spat \supseteq [\, a \ot  (\io \ot \la_N)\bigl(\al((\om_{\Ga_M(c d),x^* \Ga_M(b)} \ot
\io)(V^*) )\bigr) \mid x \in \cN, b,c \in \cN_{\psi_M}, d \in \cT, a \in A\,]
\\ & & \spat = [\, a \ot  (\io \ot \la_N)\bigl(\al((\psi_M \ot \io)(\de_M(b^* x)(c d \ot 1)) )\bigr) \mid x \in \cN, b,c \in \cN_{\psi_M}, d \in \cT, a \in A\,]
\\ & &  \spat = [\, a \ot  (\io \ot \la_N)\bigl(\al((\psi_M \ot \io)((\si^{\psi_M}_i(d) \ot 1)\de_M(b^* x) (c \ot 1)) )\bigr) \mid
\\ & & \spat\spat x \in \cN, b,c \in \cN_{\psi_M}, d \in \cT, a \in A\,] ,
\end{eqnarray*}
so that inequality (\ref{unit.ineq1}) us lets conclude that
\begin{eqnarray}
& & \hspace{-8ex} [\,(\io \ot \io \ot \la_N)((\de_M \ot \io)\al(y))\,(a \ot 1)  \mid y \in \cN_0, a \in A \, ]
\nonumber
\\ & & \hspace{-8ex} \spat \supseteq [\, a \ot  (\io \ot \la_N)\bigl(\al((\om_{p,q} \ot \io)\de_M(b^* x)
)\bigr) \mid x \in \cN, b  \in \cN_{\psi_M}, p,q \in H_M,  a \in A\,] \ . \label{unitarity.eq4}
\end{eqnarray}

\medskip

Now,
\begin{eqnarray*}
& & [\,(a^* \ot 1) \lambda v \mid a \in A, v \in H_M \ot \cK_0 \, ]
\\ & & \spat  = [\, (a^* \ot 1) \lambda \bigl(\,u \ot (\io \ot \vfi_N \ot \io)(U_{23}[\al(x_2^* x_1) \ot (\eta \ot \io)(U)])_*\,w\,\bigr) \mid
\\ & & \spat\spat a \in A, x_1,x_2 \in \cN_0, u \in H_M, w \in H_\theta\ot K, \eta \in \cA\,]
\\ & & \spat  = [\, (a^* \ot 1) (\de_M \ot \io)\bigl((\io \ot \vfi_N \ot \io)(U_{23}[\al(x_2^* x_1) \ot (\eta \ot \io)(U)])\bigr)_*\,\impl_{12}^* (u \ot w) \mid
\\ & & \spat\spat a \in A, x_1,x_2 \in \cN_0, u \in H_M, w \in H_\theta\ot K, \eta \in
\cA\,]
\\ &  & \spat = [\, (a^* \ot 1) (\de_M \ot \io)\bigl((\io \ot \vfi_N \ot \io)(U_{23}[\al(x_2^* x_1) \ot (\eta \ot \io)(U)])\bigr)_*\, (u \ot w) \mid
\\ & & \spat\spat a \in A, x_1,x_2 \in \cN_0, u \in H_M, w \in H_\theta\ot K, \eta \in
\cA\,]
\\ &  & \spat = [\, (a_2^* \ot 1) (\de_M \ot \io)\bigl((\io \ot \vfi_N \ot \io)(U_{23}[\al(x_2^* x_1) \ot (\eta \ot \io)(U)])\bigr)_*\, (a_1 u \ot w) \mid
\\ & & \spat\spat a_1,a_2 \in A, x_1,x_2 \in \cN_0, u \in H_M, w \in H_\theta\ot K, \eta \in
\cA\,]
\\ & & \spat  = [\, \bigl((a_2^* \ot 1 \ot 1) (\io \ot \io \ot \vfi_N \ot \io)\bigl(U_{34}[(\de_M \ot \io)(\al(x_2^* x_1)) \ot
(\eta \ot \io)(U)]\bigr) \\ & & \spat\spat (a_1 \ot 1 \ot 1)\bigr)_*\, (u \ot w) \mid
 a_1,a_2 \in A, x_1,x_2 \in \cN_0, u \in H_M, w \in H_\theta\ot K,
\eta \in \cA\,]
\end{eqnarray*}
For $\eta \in \cA$, we define the element $Z_\eta \in B(H_N \ot K)$ as $$Z_\eta = (T R \ot \io)(U^*) \,  (1 \ot
(\,\si^*_{\frac{i}{2}}(\bar{\eta})\, R \ot \io)(U)^*) \ .$$ By lemma \ref{integr.lem1}, we have for $a_1,a_2 \in A$,
$x_1,x_2 \in \cN_0$ and $\eta \in \cA$ that
\begin{eqnarray*}
& & (a_2^* \ot 1 \ot 1) (\io \ot \io \ot \vfi_N \ot \io)\bigl(U_{34}[(\de_M \ot \io)(\al(x_2^* x_1)) \ot (\eta \ot
\io)(U)]\bigr) (a_1 \ot 1 \ot 1) \\ & & \spat\spat\spat = (\,[(\io \ot \io \ot \la_N)((\de_M \ot \io)\al(x_2))\,(a_2
\ot1)]^*  \ot 1)
\\ & & \spat\spat \spat \spat (Z_\eta)_{34} \,(\,[(\io \ot \io \ot \la_N)((\de_M \ot \io)\al(x_1))\,(a_1 \ot 1)] \ot 1 )
\end{eqnarray*}
Similarly, using lemma \ref{integr.lem1}, we have for  $x_1,x_2 \in \cN_0$ and $\eta \in \cA$  that
\begin{eqnarray}
& & (\io \ot \vfi_N \ot \io)(U_{23}[\al(x_2^* x_1) \ot (\eta \ot \io)(U)]) \nonumber
\\ & & \spat = ((\io \ot \la_N)(\al(x_2))^* \ot 1)\,(Z_\eta)_{23} \, ((\io \ot \la_N)(\al(x_1)) \ot 1) \ . \label{unit.eq1}
\end{eqnarray}
Therefore inclusion (\ref{unitarity.eq4}) implies that
\begin{eqnarray}
& & [\,(a^* \ot 1) \lambda v \mid a \in A, v \in H_M \ot \cK_0 \, ] \nonumber
\\ & & \spat \supseteq  [\, a_2^* a_1 u \ot (\io \ot \vfi_N \ot \io)(U_{23}[\al\bigl((\om_{p_2,q_2} \ot \io)(\de_M(b_2^* x_2))^* \nonumber
\\ & & \spat \spat (\om_{p_1,q_1} \ot \io)(\de_M(b_1^* x_1))\bigr) \ot (\eta \ot \io)(U)])_* \, w  \mid a_1,a_2 \in A, x_1,x_2 \in \cN, \nonumber
\\ &  & \spat \spat  b_1,b_2 \in \cN_{\psi_M}, u \in H_M, w \in H_\theta\ot K, \eta \in \cA, p_1,p_2,q_1,q_2 \in H_M \, ]
\label{unit.eq2}
\end{eqnarray}

\smallskip

Take a bounded net $(f_i)_{i \in I}$ in $\cN_{\psi_M}$ such that $(f_i)_{i \in I}$ converges strongly$^*$ to 1.

\smallskip

Let $p,q \in H_M$ and $x \in \cN$. Then
\begin{eqnarray*}
& & (\io \ot \la_N)\bigl(\al((\om_{p,q} \ot \io)\de_M(f_i^* x))\bigr) =  (\theta_q^* \ot 1 \ot 1)(\io \ot \io \ot
\la_N)((\io \ot \al)\de_M(f_i^* x))(\theta_p \ot 1 \ot 1)
\\ & & \spat =  (\theta_q^* \ot 1 \ot 1) (\io \ot \al)(\de_M(f_i^*)) (\io \ot \io \ot
\la_N)((\io \ot \al)\de_M(x))(\theta_p \ot 1 \ot 1) \ ,
\end{eqnarray*}
so that the normality of $(\io \ot \al)\de_M$ implies that  $\bigl(\,(\io \ot \la_N)\bigl(\al((\om_{p,q} \ot
\io)\de_M(f_i^* x))\bigr)\,\bigr)_{i \in I}$ is a bounded net that converges strongly$^*$ to $(\io \ot
\la_N)\bigl(\al((\om_{p,q} \ot \io)\de_M(x))\bigr)$

\smallskip

Using this fact and formula (\ref{unit.eq1}),  we get for all $x_1,x_2 \in \cN$ and $p_1,p_2,q_1,q_2 \in H_M$ that the
net $$\bigl(\,(\io \ot \vfi_N \ot \io)(U_{23}[\al\bigl((\om_{p_2,q_2} \ot \io)(\de_M(f_i^* x_2))^*
 (\om_{p_1,q_1} \ot \io)(\de_M(f_i^* x_1))\bigr) \ot (\eta \ot \io)(U)])\,\bigr)_{i \in I}$$ is bounded and converges strongly$^*$ to
$$(\io \ot \vfi_N \ot \io)(U_{23}[\al\bigl((\om_{p_2,q_2} \ot \io)(\de_M( x_2))^*
 (\om_{p_1,q_1} \ot \io)(\de_M( x_1))\bigr) \ot (\eta \ot \io)(U)]) \ .$$

Therefore result \ref{car.res2} and inclusion (\ref{unit.eq2}) above imply that
\begin{eqnarray*}
& & [\,(a^* \ot 1) \lambda v \mid a \in A, v \in H_M \ot \cK_0 \, ]
\\ & & \spat \supseteq  [\,  u \ot (\io \ot \vfi_N \ot \io)(U_{23}[\al\bigl((\om_{p_2,q_2} \ot \io)(\de_M(x_2))^*
\\ & & \spat \spat (\om_{p_1,q_1} \ot \io)(\de_M(x_1))\bigr) \ot (\eta \ot \io)(U)])_* \, w  \mid  x_1,x_2 \in \cN,  u \in H_M,
\\ & & \spat \spat w \in H_\theta\ot K, \eta \in \cA,
p_1,p_2,q_1,q_2 \in H_M \, ]
\end{eqnarray*}

Because $M$ is supposed to be in standard form, this becomes
\begin{eqnarray*}
& & [\,(a^* \ot 1) \lambda v \mid a \in A, v \in H_M \ot \cK_0 \, ]
\\ & & \spat \supseteq  [\,  u \ot (\io \ot \vfi_N \ot \io)(U_{23}[\al\bigl((\om_2 \ot \io)(\de_M(x_2))^*
\\ & & \spat \spat (\om_1 \ot \io)(\de_M(x_1))\bigr) \ot (\eta \ot \io)(U)])_* \, w  \mid  x_1,x_2 \in \cN,  u \in H_M,
\\ & & \spat \spat w \in H_\theta\ot K, \eta \in \cA,
\om_1,\om_2 \in M_* \, ]
\\ & & \spat = [\, u \ot (\io \ot \vfi_N \ot \io)(U_{23}[\al(y_2^* y_1) \ot (\eta \ot \io)(U)])_*\,w \mid
\\ & & \spat\spat u \in H_M, w \in H_\theta \ot K, y_1,y_2 \in \cN_0, \eta \in \cA \, ]
\\ & & \spat = H_M \ot \cK_0 \ .
\end{eqnarray*}
\end{proof}

\medskip

\begin{lemma}
We have that $\lambda (H_M \ot \cK_0) = H_M \ot \cK_0$.
\end{lemma}
\begin{proof}
Define $p$ to be the projection on $\cK_0$. Also define $\phi = \lambda (1 \ot p)$, then $\phi^* \phi = (1 \ot p)
\lambda^* \lambda (1 \ot p)  = 1 \ot p$ implying that $\phi$ is a partial isometry in $M \ot B(\cK)$. Define $P =
\phi\, \phi^* =   \lambda (1 \ot p) \lambda^* \in M \ot B(\cK)$ to be the final projection of $\phi$. Thus, since
$(\de_M \ot \io)(\lambda) = \lambda_{23} \lambda_{13}$, we get that $$(\de_M \ot \io)(P) =  (\de_M \ot \io)(\lambda (1
\ot p) \lambda^*) =
 \lambda_{23} \lambda_{13} (1 \ot 1 \ot p)\lambda_{13}^* \lambda_{23}^* = \lambda_{23}  P_{13}  \lambda_{23}^* \ .$$

Because $\lambda(H_M \ot \cK_0) \subseteq H_M \ot \cK_0$, we have that $\lambda(1 \ot p) = (1 \ot p) \lambda (1 \ot
p)$, implying that $P (1 \ot p)  = \lambda (1 \ot p) \lambda^* (1 \ot p) = \lambda (1 \ot p) \lambda^* = P$. Hence
\begin{eqnarray*}
& & (\de_M \ot \io)(P) \, P_{23}  =  \lambda_{23}  P_{13}  \lambda_{23}^*  \lambda_{23} (1 \ot 1 \ot p) \lambda_{23}^*
\\ & & \spat  =  \lambda_{23}  P_{13}  (1 \ot 1 \ot p)  \lambda_{23}^* = \lambda_{23}  P_{13}  \lambda_{23}^* = (\de_M \ot \io)(P) \
.
\end{eqnarray*}
Arguing as in the proof of lemma 6.4 of \cite{VaKust}, we conclude from this that $(\de_M \ot \io)(P) = P_{23}$. So we
get for every $\om \in B(\cK)_*$ that $\de_M((\io \ot \om)(P)) = 1 \ot (\io \ot \om)(P)$ and thus $(\io \ot \om)(P) \in
\C 1$ by result 5.13 of \cite{VaKust}. This implies that $P \in (B(H_M) \ot 1)' = 1 \ot B(\cK)$. Therefore there exists
a closed subspace $\cK_1$ of $\cK$ such that $P$ is the orthogonal projection on $H_M \ot \cK_1$, thus $\lambda (H_M
\ot \cK_0)= \phi(H_M \ot \cK_0) =  H_M \ot \cK_1$. Hence the previous proposition implies that $$H_M \ot \cK_0 = [\,(a
\ot 1)\lambda v \mid a \in M, v \in H_M \ot \cK_0 \,] = [\,(a \ot 1)w \mid a \in M, w \in H_M \ot \cK_1 \,] = H_M \ot
\cK_1 ,$$ thus $\cK_1 = \cK_0$ and the lemma follows.
\end{proof}

\medskip

Now we only have to tie up the loose ends to arrive at our final conclusion (we still work under the integrability
condition of the beginning of the section!):

\begin{proposition}
The induced corepresentation $\rho$ associated to  the quadruple $(M,\de_M)$, $(N,\de_N)$, $\al$, $U$ with respect to
$(H_\theta,\pi_\theta,\la_\theta)$ is a unitary element in $M \ot B(\cK)$.
\end{proposition}
\begin{proof}
By lemma \ref{unitarity.lem1}  and the previous lemma, we know that $$\lambda(H_M \ot \cK) \subseteq H_M \ot \cK_0
\subseteq \lambda(H_M \ot \cK_0) \ .$$  Since $\lambda$ is an isometry, this implies that $H_M \ot \cK \subseteq H_M
\ot \cK_0$ and hence $\cK = \cK_0$. Using the previous lemma once more, we get that $\lambda(H_M \ot \cK) = H_M \ot
\cK$ so that $\lambda$ is unitary. Because $\rho = \lambda^*$, the proposition follows.
\end{proof}

\smallskip

It is also worthwhile remembering that $\cK  = \cK_0$.

\bigskip\medskip

\section{A correspondence between weights on $Q$ and certain weights on $M$}

\bigskip

Also in this section we will assume that $\al$ is integrable. Extend the function $\cM^+ \rightarrow Q^+ : x \mapsto
T_\al(x)$ to a function $\cT_\al : M^+ \rightarrow M^+\Ext$ such that $\cT_\al(x) = (\io \ot \vfi_N)(\al(x))$ for all
$x \in M^+$. Here we consider the map $\io \ot \vfi_N : (M \ot N)^+ \rightarrow M^+\Ext$ as an operator valued weight.
The positive extended part $Q^+\Ext$ is naturally embedded in $M^+\Ext$ (see proposition 1.9 of \cite{Haa1}). It is
proven in proposition 1.3 of \cite{Va1} that, under this embedding, $\cT_\al(M^+) \subseteq Q^+\Ext$ and that $\cT_\al
: M^+ \rightarrow Q^+\Ext$ is a semi-finite operator valued weight. We use this operator valued weight $\cT_\al$ to
pull weights down from $Q$ to $M$ (we will use proposition 2.3 of \cite{Haa1} for this):

\smallskip

\begin{definition} \label{cor.def1}
Consider a n.s.f.~weight $\eta$ on $Q$. Then we define the n.s.f.~weight $\tilde{\eta}$ on $M$ such that
$\tilde{\eta}(x) = \eta(\cT_\al(x))$ for all $x \in M^+$.
\end{definition}

For any normal weight $\eta$ on $Q$ and any element $x \in \hat{Q}_+$ the element $\eta(x) \in [0,\infty]$ is defined
in such a way that  the following holds. Let $(\eta_i)_{i \in I}$ be a family of elements  in $Q^+_*$  such that
$\eta(y) = \sum_{i \in I} \eta_i(y)$ for all $y \in Q^+$ (such a family always exist). Then $\eta(x) = \sum_{i \in I}
\eta_i(x)$.

\medskip\smallskip

Such pulled down weights satisfy a natural invariance condition with respect to $\al$ (see proposition 2.8 of
\cite{Enock}).

\smallskip

\begin{proposition}
Consider a n.s.f.~weight $\eta$ on $Q$, $a \in \cM_{\tilde{\eta}}$ and $v,w \in D(\sde_N^{\frac{1}{2}})$. Then $(\io
\ot \om_{v,w})\al(a)$ belongs to $\cM_{\tilde{\eta}}$ and $$\tilde{\eta}\bigl((\io \ot \om_{v,w})\al(a)\bigr) = \langle
\sde_N^{\frac{1}{2}} v , \sde_N^{\frac{1}{2}} w \rangle \, \tilde{\eta}(a)\ .$$
\end{proposition}

\smallskip

A similar result holds for $\psi_M$:

\begin{proposition}
Consider $a \in \cM_{\psi_M}$ and $\om \in N_*$. Then $(\io \ot \om)\al(a) \in \cM_{\psi_M}^+$ and $$\psi_M((\io \ot
\om)\al(a)) = \psi_M(a) \, \om(1) \ .$$
\end{proposition}
\begin{proof}
We may assume that $a \geq 0$ and $\om \geq 0$. We have for $\om' \in M_*^+$ that $$(\io \ot \om')\de_M((\io \ot
\om)\al(a)) = (\io \ot \om' \ot \om)((\de_M \ot \io)\al(a)) = (\io \ot (\om' \ot \om)\al)\de_M(a) \ .$$ Therefore the
right invariance of $\psi_M$ implies that $(\io \ot \om')\de_M((\io \ot \om)\al(a)) \in \cM_{\psi_M}^+$.

Using proposition 5.14 of \cite{VaKust}, we conclude from this that $(\io \ot \om)\al(a)$ belongs to $\cM_{\psi_M}^+$.
By right invariance of $\psi_M$ we have moreover for all $\om' \in M_*^+$ that
\begin{eqnarray*}
\psi_M((\io \ot \om)\al(a)) \, \om'(1) & = & \psi_M\bigl((\io \ot \om')\de_M((\io \ot \om)\al(a)))\bigr) =
\psi_M\bigl((\io \ot (\om' \ot \om)\al)\de_M(a)\bigr) \\ & = & \psi_M(a)\, (\om' \ot \om)(\al(1))=
\psi_M(a)\,\om(1)\,\om'(1) \ ,
\end{eqnarray*}
implying that $\psi_M((\io \ot \om)\al(a)) = \psi_M(a) \, \om(1)$.
\end{proof}

\medskip

Consider a n.f.s.~weight $\phi$ on $M$ and a positive self-adjoint operator $\gamma$ affiliated to $N$ such that
$\de_N(\gamma) = \gamma \ot \gamma$ and for all $a \in \cM_\phi$ and all $v,w \in D(\gamma^{\frac{1}{2}})$, we have
that the element $(\io \ot \om_{v,w})\al(a)$ belongs to $\cM_\phi$ and $$\phi\bigl((\io \ot \om_{v,w})\al(a)\bigr) =
\phi(a)\, \langle \gamma^{\frac{1}{2}} v , \gamma^{\frac{1}{2}} w \rangle \ .$$ Choose also a GNS-construction
$(H_\phi,\pi_\phi,\la_\phi)$ for $\phi$. In the next paragraphs we will state some results without proof because we
believe the reader has acquired the necessary skills and techniques by now to check these results him or herself.

\medskip

First of all, we have for $a \in \cN_\phi$, $v \in D(\gamma^{\frac{1}{2}})$ and $w \in H_N$ that $(\io \ot
\om_{v,w})\al(a) \in \cN_\phi$ and $$\|\la_\phi((\io \ot \om_{v,w})\al(a))\| \leq \|\la_\phi(a)\| \,
\|\gamma^{\frac{1}{2}} v \|\,\|w\| \ .$$

\smallskip

Let $(e_i)_{i \in I}$ be an orthonormal basis for $H_N$. Then we have for $a \in \cN_{\phi}$ and $v \in
D(\gamma^{\frac{1}{2}})$ that $$\sum_{i \in I} \|\la_\phi((\io \ot \om_{v,e_i})\al(a))\|^2 = \|\la_\phi(a)\|^2 \,
\|\gamma^{\frac{1}{2}} v\|^2 < \infty \ .$$

Therefore we can define an isometry $V_\phi \in B(H_\phi \ot H_N)$ such that $$V_\phi \, (\la_\phi(a) \ot v) = \sum_{i
\in I} \la_\phi\bigl((\io \ot \om_{\gamma^{-\frac{1}{2}} v , e_i})\al(a)\bigr) \ot e_i$$ for all $a \in \cN_\phi$ and
$v \in D(\gamma^{-\frac{1}{2}})$.

\smallskip

It follows that $(\io \ot \om_{v,w})(V_\phi)\,\la_\phi(a) = \la_\phi\bigl((\io \ot \om_{\gamma^{-\frac{1}{2}}
v,w})\al(a)\bigr)$ for all $a \in \cN_\phi$, $v \in D(\gamma^{-\frac{1}{2}})$ and $w \in H_N$.

\smallskip

Using the results of proposition 2.9 of \cite{Enock}, we  get that $V_\phi$ is a unitary element in $B(H_\phi) \ot N$
such that
\begin{enumerate}
\item $(\pi_\phi \ot \io)(\al(a)) V_\phi = V_\phi (\pi_{\phi}(a) \ot 1) $ for all $a \in M$.
\item $(\io \ot \de_N)(V_\phi) = (V_\phi)_{12} (V_\phi)_{13}$
\end{enumerate}

\bigskip

Let $\nu$ denote the scaling constant of $(N,\de_N)$ and define the strictly positive operator $P$ in $H_N$ such that
$P^{it} \la_N(a) = \nu^{\frac{t}{2}} \, \la_N(\tau^N_t(a))$ for all $t \in \R$ and $a \in \cN_{\vfi_N}$. We also let
$\hat{J}$ denote the modular conjugation of $\hat{\vfi}_N$. Recall that $\hat{J} P \hat{J} = P^{-1}$ (see proposition
2.13(2.7) of \cite{VaKust2})  and that $\tau^N_t(x) = P^{it} x P^{-it}$ and $R_N(x) = \hat{J} x^* \hat{J}$ for all $t
\in \R$ and $x \in N$.

If $\gamma$ is a strictly positive operator affiliated with $N$; arguing as in the proof of proposition 7.5 of
\cite{JK} lets us conclude that $\tau_t(\gamma) = \gamma$ for all $t \in \R$. So  $\gamma$ and $P$ strongly commute.

\smallskip

If $A$ and $B$ are strictly positive operators in $H_N$ that strongly commute, we denote by $A\!\cdot\!B$ the closure
of the composition $A B$.

\medskip

In the next proposition we will rely on relative modular theory (see e.g. section 3.11 of \cite{Stramod}).

\begin{proposition}
Let $i \in \{1,2\}$. Consider a n.f.s.~weight $\phi_i$ on $M$ with GNS-construction
$(H_{\phi_i},\pi_{\phi_i},\la_{\phi_i})$ and suppose there exists a positive self-adjoint operator $\gamma_i$
affiliated with $N$ such that $\de_N(\gamma_i) = \gamma_i \ot \gamma_i$ and for all $a \in \cM_{\phi_i}$ and all $v,w
\in D(\gamma_i^{\frac{1}{2}})$, we have that the element $(\io \ot \om_{v,w})\al(a)$ belongs to $\cM_{\phi_i}$ and
$$\phi_i\bigl((\io \ot \om_{v,w})\al(a)\bigr) = \phi_i(a)\, \langle \gamma_i^{\frac{1}{2}} v , \gamma_i^{\frac{1}{2}} w
\rangle \ .$$

Let $\nab$ denote the modular operator of $\phi_2$,$\phi_1$ with respect to the above GNS-constructions. Then
$$V_{\phi_1} \, (\nab \ot \gamma_1^{-1} \!\cdot\! P^{-1}) =   (\nab \ot \gamma_2^{-1} \!\cdot\! P^{-1})  \, V_{\phi_1}\
.$$
\end{proposition}
\begin{proof}

\smallskip

Let $T$ denote the densely defined closed linear map from within $H_{\phi_1}$ into $H_{\phi_2}$ such that
$\la_{\phi_1}(\cN_{\phi_1} \cap \cN_{\phi_2}^*)$ is a core for $T$ and  $T \la_{\phi_1}(x)  = \la_{\phi_2}(x^*)$ for
all $x \in  \cN_{\phi_1} \cap \cN_{\phi_2}^*$. Recall that the pair $J$,$\nab$ is by definition the polar decomposition
of $T$, i.e. $T = J \nab^{\frac{1}{2}}$.

\medskip

Choose $v \in D(\gamma_1^{-\frac{1}{2}})$ and $w \in D(\gamma_2^{\frac{1}{2}})$. Let $a \in \cN_{\phi_1} \cap
\cN_{\phi_2}^*$. Then $(\io \ot \om_{\gamma_1^{-\frac{1}{2}} v,w}\,)\al(a)$ belongs to $\cN_{\phi_1}$ and
$\la_{\phi_1}\bigl((\io \ot \om_{\gamma_1^{-\frac{1}{2}} v,w\,})\al(a)\bigr) = (\io \ot \om_{v,w})(V_{\phi_1})
\la_{\phi_1}(a)$. Moreover, $(\io \ot \om_{\gamma_1^{-\frac{1}{2}} v,w\,})(\al(a))^* = (\io \ot \om_{w ,
\gamma_1^{-\frac{1}{2}} v\,})(\al(a^*))$ which implies that $(\io \ot \om_{\gamma_1^{-\frac{1}{2}} v,w\,})(\al(a))$
belongs to  $\cN_{\phi_2}^*$. It follows that $(\io \ot \om_{v,w})(V_{\phi_1}) \la_{\phi_1}(a)$ belongs to $D(T)$ and
\begin{eqnarray*}
& & T\bigl((\io \ot \om_{v,w})(V_{\phi_1}) \la_{\phi_1}(a)) =  \la_{\phi_2}\bigl((\io \ot \om_{\gamma_1^{-\frac{1}{2}}
v,w\,})(\al(a))^*\bigr) = \la_{\phi_2}\bigl((\io \ot \om_{w , \gamma_1^{-\frac{1}{2}} v\,})(\al(a^*))\bigr)
\\
& & \spat = (\io \ot \om_{\gamma_2^{\frac{1}{2}} w , \gamma_1^{-\frac{1}{2}} v\,})(V_{\phi_2}) \la_{\phi_2}(a^*) = (\io
\ot \om_{\gamma_2^{\frac{1}{2}} w , \gamma_1^{-\frac{1}{2}} v\,})(V_{\phi_2})\, T \la_{\phi_1}(a) \ .
\end{eqnarray*}
Since such elements $\la_{\phi_1}(a)$ form a core for $T$, we conclude that
\begin{equation}
(\io \ot \om_{\gamma_2^{\frac{1}{2}} w , \gamma_1^{-\frac{1}{2}} v\,})(V_{\phi_2})\, T \subseteq T\,(\io \ot
\om_{v,w})(V_{\phi_1}) \ . \label{cor.eq2}
\end{equation}
Taking the adjoint of this inclusion we find that
\begin{equation}
(\io \ot \om_{w,v})(V_{\phi_1}^*)\,T^* \subseteq  T^*\,(\io \ot \om_{\gamma_1^{-\frac{1}{2}} v, \gamma_2^{\frac{1}{2}}
w\,})(V_{\phi_2}^*) \label{cor.eq3}
\end{equation}

\smallskip

Since $(\io \ot \de_N)(V_{\phi_i}) = (V_{\phi_i})_{12} \, (V_{\phi_i})_{13}$, we get for every  $\om \in M_*$, that the
element $(\om \ot \io)(V_{\phi_i})$ belongs to $D(S_N)$ and $S_N( (\om \ot \io)(V_{\phi_i})) = (\om \ot
\io)(V_{\phi_i}^*)$. Since $S_N  = \tau^N_{-\frac{i}{2}}\, R_N$, this implies easily that
\begin{equation}
(\io \ot \om_{v,w})(V_{\phi_i}^*) = (\io \ot \om_{ \hat{J} P^{\frac{1}{2}}  w , \hat{J} P^{-\frac{1}{2}}
v\,})(V_{\phi_i}) \ .  \label{cor.eq4}
\end{equation}
for all $v \in D(P^{-\frac{1}{2}})$ and $w \in D(P^{\frac{1}{2}})$.

\smallskip

So we get for $v \in D(\gamma_1^{\frac{1}{2}} P^{\frac{1}{2}} \gamma_1^{\frac{1}{2}} P^{\frac{1}{2}})$ and $w \in
D(\gamma_2^{-\frac{1}{2}} P^{-\frac{1}{2}} \gamma_2^{-\frac{1}{2}} P^{-\frac{1}{2}})$.

\begin{eqnarray*}
& & (\io \ot \om_{v,w})(V_{\phi_1})\,\nab  =  (\io \ot \om_{ \hat{J} P^{-\frac{1}{2}}  w , \hat{J} P^{\frac{1}{2}}
v\,})(V_{\phi_1}^*) \,T^* T  \subseteq  T^*\,(\io \ot \om_{\gamma_1^{-\frac{1}{2}} \hat{J} P^{\frac{1}{2}} v ,
\gamma_2^{\frac{1}{2}} \hat{J} P^{-\frac{1}{2}}  w \ ,})(V_{\phi_2}^*)\,T \\ &  & \spat = T^*\,(\io \ot \om_{\hat{J}
\gamma_1^{\frac{1}{2}} P^{\frac{1}{2}} v , \hat{J} \gamma_2^{-\frac{1}{2}}  P^{-\frac{1}{2}}  w \ ,})(V_{\phi_2}^*)\,T
= T^*\, (\io \ot \om_{\hat{J} P^{\frac{1}{2}} \hat{J} \gamma_2^{-\frac{1}{2}}  P^{-\frac{1}{2}}  w , \hat{J}
P^{-\frac{1}{2}} \hat{J} \gamma_1^{\frac{1}{2}} P^{\frac{1}{2}} v\,})(V_{\phi_2}) \, T \\ &  & \spat = T^*\, (\io \ot
\om_{P^{-\frac{1}{2}} \gamma_2^{-\frac{1}{2}}  P^{-\frac{1}{2}}  w ,
 P^{ \frac{1}{2}}  \gamma_1^{\frac{1}{2}} P^{\frac{1}{2}} v\,})(V_{\phi_2}) \, T
\subseteq T^* T\, (\io \ot \om_{\gamma_1^{\frac{1}{2}} P^{ \frac{1}{2}}  \gamma_1^{\frac{1}{2}} P^{\frac{1}{2}} v ,
\gamma_2^{-\frac{1}{2}} P^{-\frac{1}{2}} \gamma_2^{-\frac{1}{2}}  P^{-\frac{1}{2}}  w \,})(V_{\phi_1}) ,
\end{eqnarray*}
or in other words,
 $$(\io \ot \om_{v,w})(V_{\phi_1})\,\nab \subseteq \nab\,(\io \ot \om_{\gamma_1^{\frac{1}{2}} P^{
\frac{1}{2}} \gamma_1^{\frac{1}{2}} P^{\frac{1}{2}} v , \gamma_2^{-\frac{1}{2}} P^{-\frac{1}{2}}
\gamma_2^{-\frac{1}{2}} P^{-\frac{1}{2}}  w \,})(V_{\phi_1}) $$ Since $D(\gamma_1^{\frac{1}{2}} P^{\frac{1}{2}}
\gamma_1^{\frac{1}{2}} P^{\frac{1}{2}})$ is a core for $\gamma_1 \!\cdot\! P$ and $D(\gamma_2^{-\frac{1}{2}}
P^{-\frac{1}{2}} \gamma_2^{-\frac{1}{2}} P^{-\frac{1}{2}})$ is a core for $\gamma_2^{-1} \!\cdot\! P^{-1}$, it follows
easily that $$(\io \ot \om_{v,w})(V_{\phi_1})\,\nab \subseteq \nab \, (\io \ot \om_{(\gamma_1 \cdot P)\,v,
(\gamma_2^{-1} \cdot P^{-1})\,w\,})(V_{\phi_1}) $$ for all $v \in D(\gamma_1 \!\cdot\! P)$ and $w \in D(\gamma_2^{-1}
\!\cdot\! P^{-1})$. By lemma 5.9 of \cite{VaKust}, we conclude that
\begin{equation}
V_{\phi_1} \, (\nab \ot \gamma_1^{-1} \!\cdot\! P^{-1}) \subseteq   (\nab \ot \gamma_2^{-1} \!\cdot\! P^{-1})  \,
V_{\phi_1} \ . \label{cor.eq5}
\end{equation}
Taking the adjoint of this equation of this equation and multipliying it with $V_{\phi_1}$ from the left and the right,
we arrive at the other inclusion.
\end{proof}

\medskip

\begin{corollary}
Define the strongly continuous one-parameter group $\kappa$ on $N$ by setting $\kappa_t(x) =$ \newline $\sde_N^{it}\,
\tau^N_t(x) \, \sde_N^{-it}$ for all $t \in \R$ and $x \in N$. Then $$\al\, \si_t^{\tilde{\theta}} =
(\si_t^{\tilde{\theta}} \ot \kappa_{-t})\al$$ for all $t \in \R$.
\end{corollary}
\begin{proof}
Apply the previous proposition with $\phi_1 = \phi_2 = \tilde{\theta}$ and some GNS-construction
$(\tilde{H},\tilde{\pi},\tilde{\la})$ for $\tilde{\theta}$. In this case , $\gamma_1 = \gamma_2 = \sde_N$ and $\nab$ is
the modular operator for $\tilde{\theta}$ in this GNS-construction. Then we get for $x \in M$ and $t \in \R$ that
\begin{eqnarray*}
& & (\tilde{\pi} \ot \io)(\al(\si_t^{\tilde{\theta}}(x))) = V_{\tilde{\theta}} (\tilde{\pi}(\si_t^{\tilde{\theta}}(x))
\ot 1) V_{\tilde{\theta}}^*
\\ & & \spat = V_{\tilde{\theta}}
(\nab^{it} \tilde{\pi}(x) \nab^{-it}  \ot 1) V_{\tilde{\theta}}^* = V_{\tilde{\theta}} (\nab^{it} \tilde{\pi}(x)
\nab^{-it}  \ot \sde_N^{-it} P^{-it} P^{it} \sde_N^{it} ) V_{\tilde{\theta}}^*
\\ & & \spat = (\nab^{it} \ot \sde_N^{-it} P^{-it}) V_{\tilde{\theta}} ( \tilde{\pi}(x)
 \ot 1 ) V_{\tilde{\theta}}^* (\nab^{-it}  \ot  P^{it} \sde_N^{it} )
\\ & & \spat = (\nab^{it} \ot \sde_N^{-it} P^{-it}) (\tilde{\pi} \ot \io)(\al(x)) (\nab^{-it}  \ot  P^{it} \sde_N^{it} )
\\ & & \spat = (\tilde{\pi} \ot \io)\bigl(\,(\si^{\tilde{\theta}}_t \ot \kappa_{-t})\al(x)\,\bigr) \ ,
\end{eqnarray*}
and the corollary follows.
\end{proof}

\medskip

In the next two proposition we characterize all the weights on $M$ which can be pulled down from $Q$ by the procedure
in definition \ref{cor.def1}

\begin{proposition}
Consider a n.s.f.~weight $\eta$ on $Q$. Then $$\al\bigl((D \tilde{\eta},D \psi_M)_t\bigr) = (D \tilde{\eta},D \psi_M)_t
\ot \sde_N^{-it}$$ for all $t \in \R$.
\end{proposition}
\begin{proof}
Call $\nab$ the modular operator of the pair $\tilde{\eta}$, $\psi_M$. By equation (11) of section 3.11 of
\cite{Stramod} we know that $(D \tilde{\eta},D \psi_M)_t = \nab^{it} \nabp_M^{-it}$ for all $t \in \R$. Applying the
previous proposition to the pair $\tilde{\eta}$, $\psi_M$ (in which case $\gamma_1 = 1$ and $\gamma_2 = \sde_N$), we
get that $V_{\psi_M} (\nab \ot P^{-1}) = (\nab \ot \sde_N^{-1}\!\cdot\! P^{-1}) V_{\psi_M}$. If we apply the previous
proposition to the pair $\psi_M$, $\psi_M$ (in which case $\gamma_1 = \gamma_2 = 1$, we get that $V_{\psi_M} (\nabp_M
\ot P^{-1}) \subseteq (\nabp_M \ot P^{-1}) V_{\psi_M}$.

So we get for all $t \in \R$ that
\begin{eqnarray*}
\al\bigl((D \tilde{\eta},D \psi_M)_t\bigr) & = & V_{\psi_M} ((D \tilde{\eta},D \psi_M)_t \ot 1) V_{\psi_M}^* =
V_{\psi_M} (\nab^{it} \nabp_M^{-it} \ot P^{-it} P^{it}) V_{\psi_M}^*
\\ & = & (\nab^{it} \ot \sde_N^{-it} P^{-it} ) V_{\psi_M} V_{\psi_M}^* (\nabp_M^{-it} \ot P^{it})
= (\nab^{it} \ot \sde_N^{-it} P^{-it}) (\nabp_M^{-it} \ot P^{it}) \\ & = & (D \tilde{\eta},D \psi_M)_t \ot \sde_N^{-it}
\ .
\end{eqnarray*}
\end{proof}

\medskip

It is now easy to prove the converse of the previous result.

\begin{proposition}
Consider a n.s.f.~weight $\phi$ on $M$ such that $\al\bigl((D \phi, D \psi_M)_t) = (D \phi, D \psi_M)_t \ot
\sde_N^{-it}$ for all $t \in \R$. Then there exists a unique n.s.f.~weight $\eta$ on $Q$ such that $\phi =
\tilde{\eta}$.
\end{proposition}
\begin{proof}
By the previous proposition, we have for $t \in \R$ that
\begin{eqnarray*}
& & \al((D \phi, D \tilde{\theta})_t)  = \al\bigl((D \phi , D \psi_M)_t (D \psi_M, D\tilde{\theta})_t\bigr) \\ & &
\spat = ((D \phi , D \psi_M)_t \ot \sde_N^{-it}) \, ( (D \psi_M ,D\tilde{\theta})_t \ot \sde_N^{it}) = (D \phi, D
\tilde{\theta})_t \ot 1 \ ,
\end{eqnarray*}
from which we conclude that $(D \phi, D \tilde{\theta})_t$ belongs to $Q$. By theorem 4.7(1) of \cite{Haa1} we also get
that $$(D \phi, D \tilde{\theta})_{s+t}  = (D \phi, D \tilde{\theta})_s\, \si_s^{\tilde{\theta}}((D \phi, D
\tilde{\theta})_t) = (D \phi, D \tilde{\theta})_s\, \si_s^{\theta}((D \phi, D \tilde{\theta})_t) \ .$$ Therefore
theorem 5.1 of \cite{Stramod} implies the existence of a n.f.s.~weight $\eta$ on $Q$ such that $(D\eta,D\theta)_t = (D
\phi , D\tilde{\theta})_t$ for all $t \in \R$. Theorem 4.7(2) of \cite{Haa1} tells us that
$(D\tilde{\eta},D\tilde{\theta})_t = (D\eta,D\theta)_t = (D\phi,D\tilde{\theta})_t$ for all $t \in \R$. Hence,
$\tilde{\eta} = \phi$.
\end{proof}


\bigskip\medskip

\begin{thebibliography}{VD}

\bigskip

\bibitem{Abe} {\sc E. Abe},
Hopf algebras. {\it Cambridge Tracts in Mathematics}, {\bf 74} {\it Cambridge University Press, Cambridge} (1980).

\bibitem{Dix1}  {\sc  J. Dixmier},
Les alg\`ebres d'op\'erateurs dans l'espace hilbertien. Deuxi\`eme \'edition. {\it Gauthier-Villars, Paris} (1969).

\bibitem{Blatt} {\sc R.J. Blattner}, On induced representations. {\it Amer. J. Math.} {\bf 83} (1961), 79--98.

\bibitem{Enock} {\sc M. Enock},
Sous-facteurs interm\'{e}diaires et groupes quantiques mesur\'{e}s. {\it J. Operator Theory} {\bf 42} (1999), 305--330.

\bibitem{E}  {\sc  M. Enock \&  J.-M. Schwartz},
Kac Algebras and Duality of Locally Compact Groups. {\it Springer-Verlag, Berlin}  (1992).



\bibitem{Haa1} {\sc U. Haagerup},
Operator-valued weights in von Neumann algebras, I. {\it J. Funct. Anal.} {\bf 32} (1979), 175--206.

\bibitem{Haa2} {\sc U. Haagerup},
Operator-valued weights in von Neumann algebras, I. {\it J. Funct. Anal.} {\bf 32} (1979), 175--206.


\bibitem{JK}{\sc J. Kustermans \& A. Van Daele},
\cst-algebraic quantum groups arising from algebraic quantum groups. {\it Int. J. Math.} {\bf Vol. 8, No. 8} (1997),
1067--1139.  \#q-alg/9611023


\bibitem{Kust} {\sc J. Kustermans}, Regular C*-valued weights. {\it To appear in  Journal of Operator Theory} (1997)
\#funct-an/9703005

\bibitem{Kust1} {\sc J. Kustermans}, Locally compact quantum groups in the universal setting. {\it Preprint KU Leuven}
(2000)

\bibitem{VaKust} {\sc J. Kustermans \& S. Vaes}, Locally compact quantum groups. {\it To appear in  Ann. scient. \'{E}c. Norm. Sup.} (1999)

\bibitem{VaKust2} {\sc J. Kustermans \& S. Vaes}, Locally compact quantum groups in the von Neumann algebraic setting. {\it Preprint KU Leuven.}
(2000)\ \#math/0005219

\bibitem{VaKust1} {\sc J. Kustermans \& S. Vaes}, Weight theory for C*-algebraic quantum groups. {\it Preprint University College Cork \& KU Leuven} (1999). \#math/9901063

\bibitem{Mack1} {\sc  G.W. Mackey}, Imprimitivity for representations of locally compact groups.
{\it Proc. Nat. Ac. Sci. USA} {\bf 35} (1949), 537--545.


\bibitem{Mack2} {\sc  G.W. Mackey}, Induced representations of locally compact groups, I.
{\it Ann. Math.} {\bf 55} (1952), 101--140.

\bibitem{Mack3} {\sc  G.W. Mackey}, Induced representations of locally compact groups, II.
{\it Ann. Math.} {\bf 56} (1953), 193--221.






\bibitem{Stramod} {\sc S. Stratila}, Modular Theory in Operator Algebras. {\it Abacus Press, Tunbridge Wells, England} (1981).

\bibitem{Va2} {\sc S. Vaes}, A Radon-Nikodym theorem for von Neumann algebras. {\it To appear in J. Operator Theory} (1998).
\#math/9811122


\bibitem{Va1} {\sc S. Vaes}, The unitary implementation of a locally quantum group action. {\it Preprint KU Leuven} (2000).



\bibitem{VaVan} {\sc S. Vaes \& A. Van Daele}, Hopf \cst-algebras. {\it To appear in Proc. London Math. Soc.} (1999).
\#math/9907030


\bibitem{Wor5}  {\sc  S.L. Woronowicz}, From multiplicative unitaries to quantum groups.
{\it Int. J. Math.} {\bf Vol.~7, No.~1} (1996), 127--149.


\end{thebibliography}
\end{document}